\documentclass[11pt]{elsarticle}

\usepackage{graphicx}

\usepackage{a4wide}
\usepackage{amssymb}
\usepackage{amsmath}

\usepackage{mathrsfs}
\usepackage{bm}
\usepackage{accents} 
\usepackage[latin1]{inputenc}
\usepackage[english]{babel}

\usepackage{multirow}
\usepackage{booktabs}

\usepackage{macros-rohan}

\usepackage{appendix}

\usepackage[numbers]{natbib}

\usepackage{color}

\newtheorem{theorem}{Theorem}[section]
\newtheorem{lemma}[theorem]{Lemma}
\newtheorem{e-proposition}[theorem]{Proposition}

\newtheorem{e-definition}[theorem]{Definition\rm}

\def\pathContactR{figs-R} 
\def\pathContactH{figs-H}

\newcommand{\REMOVE}[1]{} 

\newcommand\blue[1]{{#1}}

\newcommand\chER[1]{\blue{#1}}

\newcommand\chEM[1]{\blue{#1}}

\newcommand\RELEASE[2]{}

\newcommand\MeanY[1]{\Mcal_{Y}{\left (#1\right )}}

\newcommand\Hpdbav{\wtilde\Hbm_{\#}^1}

\newcommand\Proj[2]{{\rm Proj}_{#2}\left({#1}\right)}
\newcommand\Ccone{{\rm{C}}}
\newcommand\Kset{{\rm{K}}}

\def\unit#1{\,\mathrm{#1}}

\def\aOmeps#1#2{a_{\Omega}^\veps\!\left({#1},{#2}\right)}

\def\ipYs#1#2{\left \langle{#1},\,{#2}\right \rangle_{Y_s}}

\def\ipGc#1#2{\left \langle{#1},\,{#2}\right \rangle_{\Gamma_c}}
\def\iprod#1#2{\left \langle{#1}\Big |\,{#2}\right \rangle}

\def\Bllangle{\Big\langle \kern-0.35em \Big\langle}
\def\Brrangle{\Big\rangle \kern-0.35em \Big\rangle}

\def\ipprod#1#2{\Bllangle{#1}\Big |\,{#2}\Brrangle}

\def\aYs#1#2{a_{Y_S} \left ({#1},\,{#2}\right )}

\newcommand\mic{{\rm mic}}

\newcommand\nrm[2]{\left\| {#1}\right\|_{#2}}  

\newcommand\eeb[1]{\eb({#1})}

\newcommand\zerobm{\bf{0}}

\newcommand\dlt{\delta}

\newcommand\rhs{{\rm right hand side~}}
\newcommand\lhs{{\rm left hand side~}}
\newcommand\wrt{{\rm with respect to~}}
\newcommand\ie{{\textit{i.e.~}}}

\newcommand\dist{\rm{dist}}

\newcommand\jump[1]{\left[{#1}\right]_n^\veps}
\newcommand\jumpp[1]{\Big|\kern-0.16em \Big[{#1}\Big]\kern-0.16em\Big|}
\newcommand\jumpY[1]{\left[{#1}\right]_n^Y}

\newcommand\jumpYt[2]{\left[{#1}\right]_n^{({#2})}}
\newcommand\jumpYi[2]{\left[{#1}\right]_n^{{#2}}}

\def\ENorm#1#2{|\kern-0.14em|\kern-0.14em |{#1} |\kern-0.14em |\kern-0.14em|_{#2}} 

\numberwithin{equation}{section}
\def\intY{\:\sim \kern-1.17em \int}
\def\intYs{\:- \kern-0.95em \int}

\newcounter{definition}
{\vspace{7pt} \noindent{\bf Definition
  \refstepcounter{definition}\thedefinition. \protect\label{#1} \hspace{-2mm}}\rm }%
{\vspace*{-7pt} \flushright $\triangle$\\ } 

\newenvironment{myremark}[1]%
{\vspace{7pt} \noindent{\bf Remark
  \refstepcounter{remark}\theremark. \protect\label{#1} \hspace{-2mm}}\rm }%
{\vspace*{-7pt} \flushright $\triangle$\\ } 

\newcounter{proposition}
\newenvironment{myproposition}[1]%
{\vspace{7pt} \noindent{\bf Proposition
  \refstepcounter{proposition}\theproposition. \protect\label{#1} \hspace{-2mm}}\rm }%
{\vspace*{-7pt} \flushright $\triangle$\\ } 

\newenvironment{myproof}[1]%
{\vspace{7pt} \noindent{\bf Proof} of {#1} \hspace{-2mm}\rm }%
{\vspace*{-7pt} \flushright $\square$\\ } 

\begin{document}

\selectlanguage{english}

\begin{frontmatter}
%
\title{Homogenization and numerical algorithms for two-scale modelling of porous media with self-contact in micropores}


\author[1]{Eduard Rohan}
\ead{rohan@kme.zcu.cz}
\author[1]{Jan Heczko}
\ead{jheczko@ntis.zcu.cz}

\address[1]{European Centre of Excellence, NTIS -- New Technologies for
Information Society, Faculty of Applied Sciences, University of West
Bohemia, Univerzitn\'{\i} 8, 30614 Pilsen, Czech Republic}


\begin{abstract}

  The paper presents two-scale numerical algorithms for stress-strain analysis of porous media featured by self-contact at pore level.  The porosity is constituted as a periodic lattice generated by a representative cell consisting of elastic skeleton, rigid inclusion and a void pore. Unilateral frictionless contact is considered between opposing surfaces of the pore. For the homogenized model derived in our previous work, we justify incremental formulations and propose several variants of two-scale algorithms which commute iteratively solving of   the micro- and the macro-level contact subproblems. A dual formulation which take advantage of the assumed microstructure periodicity and a small deformation framework, is derived for the  contact problems at the micro-level. This enables to apply the semi-smooth Newton method. For the global, macrolevel step two alternatives are tested; one relying on a frozen contact identified at the microlevel, the other based on a reduced contact associated with boundaries of contact sets. Numerical examples of 2D deforming structures are presented as a proof of the concept. 
  
\end{abstract}

\begin{keyword}
Unilateral contact \sep Homogenization \sep Porous media \sep Variational inequality \sep Dual formulation \sep two-scale iterative algorithm 
\end{keyword}

\end{frontmatter}

\section{Introduction}\label{sec-intro}
Although the unilateral interaction between compliant bodies belongs to classic topics in structural mechanics  \cite{Hlavacek-HNL-1988solution,Haslinger-Hlavacek-Necas-HNA1996}, and efficient numerical methods have been developed, in the context of porous media and multiscale modelling, the  self-contact at pore level of such periodically heterogeneous structures present still a rather challenging problem.
Variational formulation of the unilateral contact problem in the periodic, or quasi-periodic media was treated by the homogenization techniques in a number of works \cite{Mikelic-contact1998,Cioranescu-etal-contactAA2013,Griso-etal-contact-layer-AA2016}. In our previous paper \cite{Rohan-Heczko-CaS} devoted to this subject we proposed a two-scale algorithm which is based alternating micro- and macro-level steps. There, using the asymptotic analysis approach to the homogenization, the limit problem of the unilateral frictionless self-contact in pores of the elastic skeleton has been derived, consisting of two parts. The local problems defined for any macroscopic position within the macroscopic body are formulated in terms of the variational inequality. The local responses are driven by the macroscopic strain. Based on the local true (active) contact boundary, the consistent tangent stiffness can be determined locally by solving a linear problem with the bilateral sliding contact constraint.
Consequently the global (macroscopic) problem involving the tangent stiffness tensors is solved for the macroscopic displacement field increments

In the present paper, we propose and test new modifications of the original two-scale computational algorithm reported in \cite{Rohan-Heczko-CaS}. As a novelty, a dual formulation of the pore-level contact problems in the local representative cells provides actual active contact sets which enables to compute consistent effective elastic coefficients at particular macroscopic points. At the macroscopic level, a sequential linearization leads to an incremental equilibrium problem which is constrained by a projection arising from the homogenized contact constraint, such that the Uzawa algorithm can be used. Several modifications are proposed which are related to the restriction of the assumed contact set variation in the context of the tangential incremental modulus computation.
At the local level, the finite element discretized contact problem attains the form of a nonsmooth equation which which is solved using the semi-smooth Newton method \cite{DeLuca-Facchinei-Kanzow-MathProg1996} without any regularization, or a problem relaxation. Numerical examples of 2D deforming structures are presented.

The paper is organized, as follows. In Section~\ref{sec-model}, the micromodel of elastic periodic porous structure with the unilateral frictionless contact constrain is introduced. The limit two-scale problem is recorded in Section~\ref{sec-two-scale} where the consistency of the incremental formulation is stated by virtue of Proposition~\ref{prop1}. In Section~\ref{sec-ma-mcp}, we introduce the macroscopic contact method (MCM) as the modification of the two-scale solution algorithm proposed in  \cite{Rohan-Heczko-CaS}. Some variants of the MCM algorithm are defined. The dual formulation of the local microscopic contact problems is established in Section~\ref{sec-micro}. For illustration of the proposed algorithms, numerical examples are reported in Section~\ref{sec-numex}. 

\paragraph{Notation and functional spaces}
\label{sec-Appendix-A} In the paper, the following notations are used.
\begin{list}{$\bullet$}{}




\item By $\pd_i = \pd_i^x$ we abbreviate the partial derivative $\pd/\pd x_i$.
    We use $\nabla_x = (\pd_i^x)$ and $\nabla_y = (\pd_i^y)$ when differentiation \wrt coordinate $x$ and $y$ is used, respectively.
    The symmetric gradient of a vector function
    $\ub$, $\eeb{\ub} = 1/2[(\nabla\ub)^T + \nabla\ub]$ where the transpose
    operator is indicated by the superscript ${}^T$.


\item The Lebesgue spaces $L^2(D)$ of square-integrable functions on $D$. The
    Sobolev space $H^1(D)$ of the square-integrable functions up to the 1st order generalized derivative. The notation with bold and
    non-bold letters, \ie like $H^1(D)$ and $\Hdb(D)$ is used to distinguish
    between spaces of scalar, or vector-valued functions. 

\item The space $\Hpdb(D)$ is the Sobolev space of functions defined in $D$,
    integrable up to the 1st order of the generalized derivative, and which are
    Y-periodic in $\Om$. 



    \item $\SS_2$ symmetric 2nd order tensors, $\tau_{ij} = \tau_{ji}$.
\end{list}{}{}

\section{Micro-model}\label{sec-model}
We consider porous elastic media constituted as periodic structures which can be represented by the so called representative periodic cells (RPC). In such a RPC, the pore geometry admits the unilateral self-contact while deforming the global structure. 
In this section we introduce a micromodel describing deformation of these kind of structures whose the microstructure is characterized by the scale parameter $\veps = \ell_\mic/L$, where $\ell_\mic$ and $L$ are the characteristic lengths of the microstructure and the macroscopic body.

\subsection{Porous structure and periodic geometry}\label{sec-geom}
An open bounded domain $\Om\subset \R^d$, with the dimension $d=2,3$, is constituted by the solid skeleton $\Om_s$ and by the fractures (fissures) $\Om_f$, so that
\begin{equation}\label{eq-g1}
\begin{split}
  \Om & = \Om_s^\veps \cup \Om_f^\veps \cup \Gamma^\veps\;, \quad \Om_s^\veps \cap \Om_f^\veps = \emptyset\;,\quad \ol{\Om_f^\veps} \subset \Om\;,
\end{split}
\end{equation}
where $\Gamma^\veps = \ol{\Om_s^\veps} \cap \ol{\Om_f^\veps}$ is the interface.
Further we assume that $\Om_s^\veps$ is a connected domain, whereas $\Om_f^\veps$ may not be connected.
To impose boundary conditions, the decomposition of the boundary is introduced, as follows:
\begin{equation}\label{eq-g2}
\begin{split}
  \pd\Om & = \pd_u\Om \cup \pd_\sigma \Om\;,\quad \pd_u\Om \cap \pd_\sigma \Om = \emptyset\;,\\
  \Gamma_c^\veps & = \Gamma^\veps \setminus \pd_\sigma \Om_s^\veps\;,\\
  \Gamma_+^\veps & =  \Gamma_c^\veps \setminus \Gamma_-^\veps\;,\quad \Gamma_+^\veps \cap  \Gamma_-^\veps = \emptyset\;,
\end{split}
\end{equation}
where $\pd_\sigma \Om_s^\veps \cap \Gamma^\veps$ is a part of the interface $\Gamma^\veps$ on which any contact is excluded; note that $\pd_\sigma \Om \subset \pd_\sigma \Om_s^\veps$.
Boundary $\Gamma_c^\veps$ splits into two disjoint parts, such that, in the deformed configuration,  points on $\Gamma_+^\veps$ can get in contact with those situated on $\Gamma_-^\veps$. 

	

	

\begin{figure}[h]
	\centering
\includegraphics[width=0.8\linewidth]{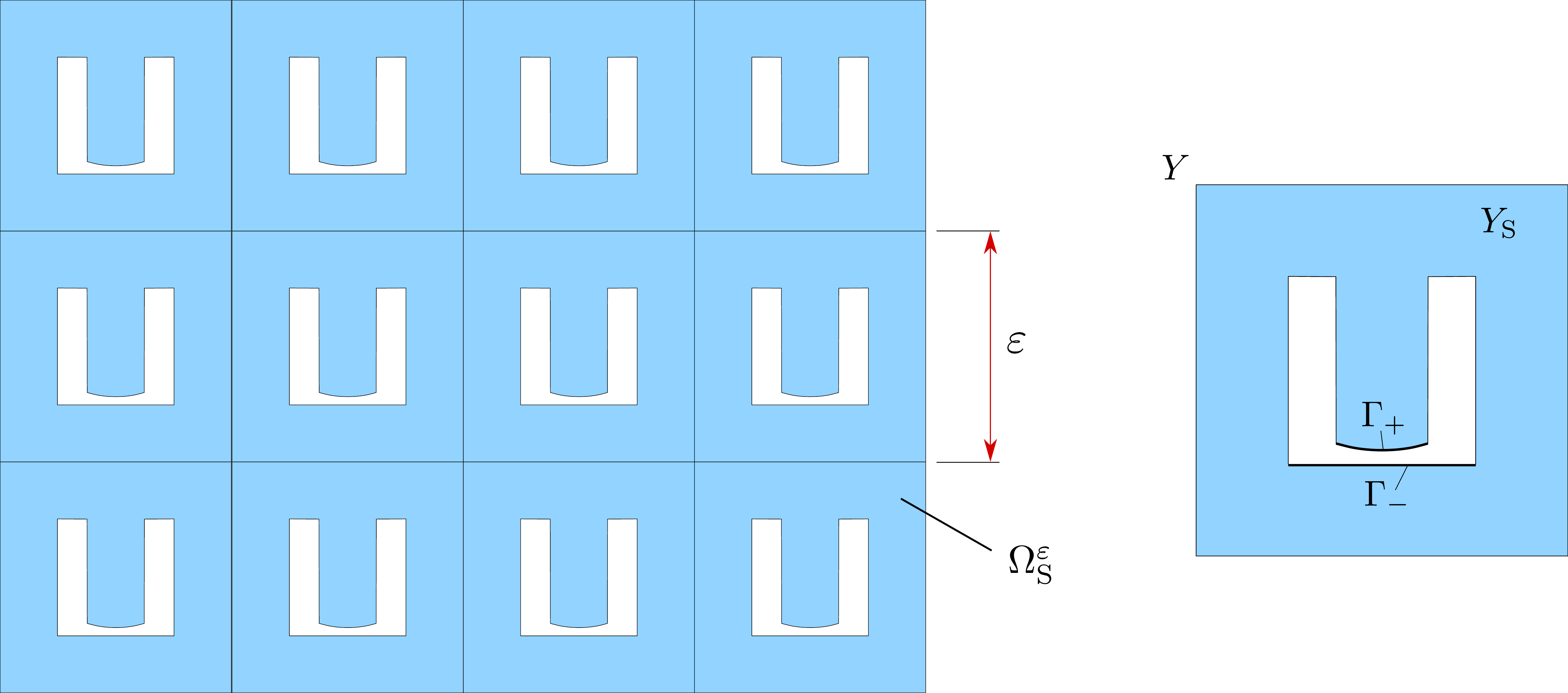} 

\caption{Periodic structure with rigid inclusions. Right: representative periodic cell $Y$ }\label{fig-OmY}
	
\end{figure}

The solid part $\Om_s^\veps$ is generated as a periodic lattice by repeating the representative volume element (RVE) occupying domain $\veps Y$. The zoomed cell 
$Y = \Pi_{i=1}^3]0,\bar y_i[ \subset \RR^3$ 
splits into the solid part
occupying domain $Y_s$ and the complementary fissure part $Y_f$, see Fig.~\ref{fig-OmY}, thus
\begin{equation}\label{eq-mi6}
\begin{split}
Y   = Y_s  \cup Y_f  \cup \Gamma^Y \;,\quad
Y_s   = Y \setminus Y_f  \;,\quad
\Gamma_Y   = \ol{Y_s } \cap \ol{Y_f }\;.
\end{split}
\end{equation}
For a given scale $\veps > 0$, $\ell_i = \veps \bar y_i$ is the characteristic size associated with the $i$-th coordinate direction, whereby
also $\veps \approx \ell_i / L$, hence $\ell_i \approx \ell^\mic$ (for all $i=1,2,3$) specifies the microscopic characteristic length $\ell^\mic$.
The contact boundary is subject to the analogous split as in \eq{eq-g2},
\begin{equation}\label{eq-mi6a}
\begin{split}
  \Gamma_c^Y & \subset \Gamma^Y\;,\\
  \Gamma_+^Y & =  \Gamma_c^Y \setminus \Gamma_-^Y\;,\quad \Gamma_+^Y \cap  \Gamma_-^Y = \emptyset\;.
\end{split}
\end{equation}

\subsection{Contact problem formulation}\label{sec-problem}
To introduce the contact kinematic conditions,
with reference to the contact boundary split \eq{eq-g2} and denoting by  $\nb(x^-)$ the unit normal to $\Gamma_c^\veps$ at $x^- \in \Gamma_-^\veps$, let $x^+ = \xi \nb(x^-) + x^- \in \Gamma_+^\veps$ for some $\xi\geq 0$. Thus, two matching points on the contact surfaces $\Gamma_+^\veps$ and $\Gamma_-^\veps$ are introduced, which enables to define the jump 
$\jump{~}$ and, thereby,  the contact gap function $g_c^\veps$, as follows,
\begin{equation}\label{eq-p0}
  \begin{split}
    g_c^\veps(\ub^\veps)& =  \jump{\ub^\veps} - s^\veps \;,\\
  \mbox{ where }  s^\veps & = \jump{x}\;,\\
    \jump{\ub}& = \nb(x^-)\cdot(\ub(x^+) - \ub(x^-))\;,\quad x^+ \in \Gamma_+^\veps,\quad  x^- \in \Gamma_-^\veps\;.
  \end{split}
\end{equation}


In this paper, we shall assume vanishing traction forces on the non-contact surface of the fissures, thus, $\bb^\veps \equiv 0 $ on $\pd_\sigma\Om_s^\veps\setminus \pd_\sigma\Om$. Moreover, let
 $\bb^\veps = \bb$ be independent of $\veps$ on $\pd_\sigma\Om$.

We shall now introduce the weak formulation of the friction-less contact problem for linear elastic structures subject to small strains and linearized contact conditions.
For this, the set of kinematically admissible displacements is needed,
\begin{equation}\label{eq-K}
\begin{split}
\Kcal^\veps = \{\vb \in \Hdb(\Om_s^\veps)|\;\vb = 0 \mbox{ on } \pd_u\Om_s^\veps\;,\; g_c^\veps(\vb) \leq 0\; \mbox{ on } \Gamma_c^\veps\}\;.
\end{split}
\end{equation}

\paragraph{Weak solutions of the contact problem} 
Given body forces $\fb^\veps$ and external surface traction forces $\bb$, find
a displacement field  $\ub^\veps\in \Kcal^\veps$ which satisfies the variational inequality,
\begin{equation}\label{eq-p2}
\begin{split}
 \aOmeps{\ub^\veps}{\vb^\veps - \ub^\veps} & \geq
  \int_{\Om_s^\veps}\fb^\veps\cdot(\vb^\veps - \ub^\veps) + \int_{\pd_\sigma\Om}\bb\cdot(\vb^\veps - \ub^\veps)
  \;,\quad \forall \vb^\veps \in \Kcal^\veps\;,\\
  \mbox{ where }\quad \aOmeps{\wb}{\vb} & = \int_{\Om_s^\veps} \Dop\eeb{\wb}:\eeb{\vb}\;.
\end{split}
\end{equation}

\section{Two scale problem}\label{sec-two-scale}
The homogenization of the porous medium governed by the contact problem \eq{eq-p2} has been reported in \cite{Rohan-Heczko-CaS}, where two scale problem was derived using the unfolding method using the theoretical results proved by Cioranescu \etal, \cite{Cioranescu-etal-contactAA2013}. For the sake of brevity, we shall disregard loading by surface tractions, thus, we put $\bb\equiv\zerobm$.

Since the main purpose of the present paper is to develop an effective algorithm for computing numerical solutions of the contact problem for $\veps\rightarrow 0$, we present only the homogenized model arising from formulation  \eq{eq-p2}.

where $\ub^0 \in \Hdb(\Om)$ and $\ub^1 \in L^2(\Om;\Hpdb(Y_s))$, and $\what y = y - \MeanY{y}$ is the relative position \wrt the barrycenter of $Y$. Then it is straightforward to introduce the asymptotic expansion of solutions,
\begin{equation}\label{eq-p3}
\begin{split}
  \Tuf{\ub^\veps(x)} & = \ub^0(x) + \veps \ub^1(x,y) + \veps^2(\dots\;.
\end{split}
\end{equation}
In analogy, we consider the truncated expansions of the test functions $\Tuf{\vb^\veps(x)}  = \vb^0(x) + \veps \vb^1(x,y)$, where
$\vb^0 \in \Hdb(\Om)$ and $\vb^1 \in L^2(\Om;\Hpdb(Y_s))$.

Consequently, we can introduce the kinematic admissibility set $W_\#(Y_s,\Om)$ defined in terms of the gap function $g_c^Y$ and the associated convex set
$\Kcal_Y$,
\begin{equation}\label{eq-p3b}
  \begin{split}
    W_\#(Y_s,\Om) & = \{(\vb^0,\vb^1) \in \Hdb(\Om) \times L^2(\Om;\Hpdb(Y_s))\,|\;
  g_c^Y(\vb^1,\nabla_x\vb^0) \leq 0\}\;,\\
     \Kcal_Y(\nabla \ub) & =\{\vb \in \Hpdb(Y_s)| \; g_c^Y(\vb,\nabla \ub) \leq 0\}\;, \\
  \mbox{ where }\quad g_c^Y(\ub^1,\nabla \ub^0) & = \jumpY{ \nabla \ub^0\what y + \ub^1 - \what y} \;.
\end{split}
\end{equation}

Solutions $\ub^\veps$ of problem \eq{eq-p2} converge 
to the two-scale solutions $\Scal := (\ub^0,\ub^1)$ of the following limit problem: 

For given $\bar\fb \in \Lb^2(\Om)$ 
, find o two-scale solution $(\ub^0,\ub^1)\in W_\#(Y_s,\Om)$, such that
\begin{equation}\label{eq-if1}
\begin{split}
  \int_\Om\intY_{Y_s} \Dop(\eeby{\ub^1} + \eebx{\ub^0}):(\eeby{\vb^1} + \eebx{\vb^0} - \eeby{\ub^1} - \eebx{\ub^0})
  \geq & \int_\Om \bar \fb \cdot (\vb^0 - \ub^0)\;,
\end{split}
\end{equation}
for all $(\vb^0,\vb^1) \in W_\#(Y_s,\Om)$.

We shall need the elastic bilinear form
\begin{equation}\label{eq-a}
\begin{split}
 \aYs{\ub}{\vb} & = \intY_{Y_s}\Dop\eeby{\ub}:\eeby{\vb}\;.
\end{split}
\end{equation}

In \cite{Rohan-Heczko-CaS}, we have shown that 
the two-scale solutions $(\ub^0,\ub^1)$ satisfy the decomposed limit problem \eq{eq-p4}-\eq{eq-p5}, consisting of the Local and the Global equilibria.
\begin{list}{}{}
\item Local equilibrium: for a.a. $x \in \Om$, the fluctuating displacement fields $\ub^1(x,\cdot) \in \Kcal_Y(\nabla\ub^0)$, 
satisfy the Local variational inequality (LVI), 
\begin{equation}\label{eq-p4}
\begin{split}
  \aYs{\ub^1 + \Pibf^{ij}e_{ij}^x(\ub^0)}{\vb - \ub^1}
  \geq 0\;,\quad \forall \vb \in \Kcal_Y(\nabla\ub^0)\;.
\end{split}
\end{equation}

\item Global equilibrium: Macroscopic displacement $\ub^0 \in U_0(\Om)$, where 
$U_0(\Om) = \{\vb \in \Hdb(\Om)|\;\vb = 0 \mbox{ on } \pd_u\Om\}$, 
 satisfies the Global variational equality (GVE),
  \begin{equation}\label{eq-p5}
\begin{split}
  \int_\Om \sigmabf^0(\ub^0,\ub^1): \eebx{\vb^0} & = \int_\Om \bar\fb\cdot \vb^0 + \int_{\pd_\sigma\Om} \bb\cdot \vb^0 \quad \forall \vb \in U_0(\Om)\;,\\
  \mbox{ with } \quad  \sigmabf^0 = (\sigma_{ij}^0), \quad  \sigma_{ij}^0 & =\aYs{\ub^1 + \Pibf^{kl}e_{kl}^x(\ub^0)}{\Pibf^{ij}}\;.
\end{split}
  \end{equation}

  \subsubsection{Two-scale incremental formulation}
  To solve this two-scale nonlinear problem, we propose an incremental formulation. We establish  the solution increments,
\begin{equation}\label{eq-if5a}
\begin{split}
  \sigmabf^0 & = \tilde{\sigmabf}^0 + \dlt\sigmabf^0\;,\\
  \ub^0 & = \tilde{\ub}^0 + \dlt\ub^0\;,\\
  \ub^1 & = \tilde{\ub}^1 + \dlt\ub^1\;. 
\end{split}
\end{equation}
Using the space 
\begin{equation}\label{eq-if16}
 V_0(\Ub,Y_s,x) = \{\vb \in \Hpdb(Y_s)|\; \jumpYt{\vb + \Ub}{c} = 0 \mbox{ on }\Gamma_c^*(x)\}\;,
\end{equation}
we introduce a local corrector problem: Find
$\hat\wb^{ij}(x,\cdot) \in V_0(\bmi{0},Y_s,x)$, such that
\begin{equation}\label{eq-if17}
  \begin{split}
    \aYs{\hat\wb^{ij} + \Pibf^{ij}}{\vb} & = 0\;,\\
 \end{split}    
\end{equation} 
for all $\vb  \in V_0(\bmi{0},Y_s,x)$. 
Upon introducing the tangent stiffness tensor
\begin{equation}\label{eq-if19}
  \begin{split}
 D_{ijkl}^H & = \aYs{\hat\wb^{kl}+\Pibf^{kl}}{\hat\wb^{ij}+\Pibf^{ij}}\;,
\end{split}    
\end{equation}
we define the macroscopic incremental equilibrium, 
\begin{equation}\label{eq-if20}
  \begin{split}
    \dlt{\ub}^0\in U_0(\Om)\;:\quad
 \int_\Om   D_{ijkl}^He_{kl}^x(\dlt{\ub}^0)e_{ij}^x(\tilde\vb^0)& = \int_\Om \dlt\fb\cdot\tilde\vb^0\;,\quad \forall \tilde\vb^0\in U_0(\Om)\;.
\end{split}    
\end{equation}

\begin{myproposition}{prop1}
  Let us assume the pair $(\tilde{\ub}^0,\tilde{\ub}^1)$ satisfies the LVI (local variational inequality) \eq{eq-p4} which provide a true contact set $\tilde\Gamma_c^*(x)$ for a.a. $x \in \Om$, 
 and the GVE (global variational equality)   \eq{eq-p5}.  Consider $\wtilde{\Kcal_Y}(\nabla({\ub}^0))  := {\Kcal_Y}(\nabla({\ub}^0)) - \tilde{\ub}^1$.
  \begin{list}{}{}
  \item (i) 
      An incremental formulation is consistent with limit decomposed problem \eq{eq-p4}-\eq{eq-p5} and with the limit variational inequality \eq{eq-if1}, if for a given load increment $\dlt\fb \in L^2(\Om)$ the increments $\dlt{\ub}^0 \in U_0(\Om)$ and $\dlt{\ub}^1 \in  \wtilde{\Kcal_Y}(\nabla({\ub}^0) $ verify: 
\begin{equation}\label{eq-if12}
  \begin{split}
    \aYs{\dlt\ub^1 + \Pibf^{ij}e_{ij}^x(\dlt{\ub}^0)}{\dlt\tilde\vb - \dlt\ub^1} & \geq 0\;,\quad \forall \dlt\tilde\vb \in \wtilde{\Kcal_Y}(\nabla({\ub}^0))\\
     \int_\Om\intY_{Y_s} \Dop(\eeby{\dlt\ub^1} + \eebx{\dlt\ub^0}):(\eebx{\tilde\vb^0} - \eeby{\wb(\tilde\vb^0)}) & = \int_\Om \dlt\fb\cdot\tilde\vb^0\;,\quad \forall \tilde\vb^0 \in  U_0(\Om)\;.
\end{split}    
\end{equation}
      
\item (ii) Assuming the true contact boundaries do not change with the increments, \ie the contact
$g_c^Y(\ub^1,\ub^0) = 0$ and $g_c^Y(\tilde\ub^1,\tilde\ub^0) = 0$ hold only 
  on $\tilde\Gamma_c^*(x)$ for a.a. $x \in \Om$, then the increment $(\dlt\ub^1,\dlt\ub^0)$ satisfies \eq{eq-if20} with the tangent stiffness \eq{eq-if19}.

  \end{list}
\end{myproposition}
In the proof, we consider the following substitutions 
\begin{equation}\label{eq-if3}
\begin{split}
  \vb^1 & := \tilde{\vb}^1 - \wb(\tilde{\vb}^0)\;,\\
  \vb^0 & := \tilde{\vb}^0 + \ub^0\;,\\
  \wb(\tilde{\vb}^0) & = \nabla_x \tilde{\vb}^0 \hat y \quad \mbox{ on } \Gamma_c\;,
\end{split}
\end{equation}
where $\vb^1(x,\cdot) \in \Hpdb(Y_s)$ and also $\wb(\vb)(x,\cdot) \in \Hpdb(Y_s)$ for any
$\vb \in \Hdb(\Om)$.
The following lemma holds.
\begin{lemma}\label{lem1}
Let $\ub \in \Kcal_Y(\Gb^0)$ for a given $\Gb = (G_{ij})$ and $\wb(\bar\ub) \in \Hpdb(Y_s)$ is given by \eq{eq-if3}$_3$ for any $\bar\ub \in \Hdb(\Om)$.
If $\jumpYi{\vb}{c} = 0$, then $\ub + \vb \in \Kcal_Y(\Gb)$. Furthermore,
$ \Kcal_Y(\nabla(\bar\ub) + \Gb) =  \Kcal_Y(\Gb) - \wb(\bar\ub)$.
\end{lemma}

\begin{myproof}{ Proposition~\ref{prop1}} \\
  (i)
Using \eq{eq-if3} substituted in  \eq{eq-if1}, we get
\begin{equation}\label{eq-if5}
\begin{split}
  \int_\Om \aYs{\tilde{\ub}^1 + \Pibf^{ij}e_{ij}^x(\tilde{\ub}^0)}{\tilde{\vb}^1-\wb(\tilde{\vb}^0) - \tilde{\ub}^1 - \dlt\ub^1} + \int_\Om \tilde{\sigmabf}^0(\tilde{\ub}^0,\tilde{\ub}^1):\eebx{\tilde\vb^0} \\
+  \int_\Om\intY_{Y_s} \Dop(\eeby{\dlt\ub^1} + \eebx{\dlt\ub^0}):(\eebx{\tilde\vb^0} - \eeby{\wb(\tilde\vb^0)})\\ + \int_\Om \aYs{\dlt\ub^1 + \Pibf^{ij}e_{ij}^x(\dlt{\ub}^0)}{\tilde{\vb}^1 -\tilde{\ub}^1 - \dlt\ub^1} \geq  \int_\Om (\tilde\fb + \dlt\fb)\cdot\tilde{\vb}^0\;,
\end{split}
\end{equation}
By the consequence of Lemma~\ref{lem1},
\begin{equation}\label{eq-if6}
\begin{split}
  \Kcal_Y(\nabla(\tilde{\ub}^0 + \dlt\ub^0)) &= \Kcal_Y(\nabla(\dlt\ub^0)) - \wb(\tilde{\ub}^0)
  = \Kcal_Y(\nabla\tilde\ub^0) - \wb(\dlt{\ub}^0)\;.
\end{split}
\end{equation}
Since $\tilde{\ub}^1 \in \Kcal_Y(\nabla(\tilde{\ub}^0))$ solves the LVI and also
$\tilde{\ub}^1 + \dlt\ub^1 \in \Kcal_Y(\nabla({\ub}^0))$, we can define $\tilde\wb:= \dlt\ub^1 + \wb(\dlt{\ub}^0)$ and get
\begin{equation}\label{eq-if6a}
  \begin{split}
    \tilde{\vb}^1 - \wb(\tilde{\vb}^0 + \wb(\dlt{\ub}^0) \in \Kcal_Y(\nabla\tilde\ub^0)\;,\\
   \tilde{\ub}^1 + \tilde\wb = \tilde{\ub}^1 + \dlt\ub^1 + \wb(\dlt{\ub}^0)  \in \Kcal_Y(\nabla\tilde\ub^0)\;.
\end{split}
\end{equation}
Due to \eq{eq-if6a}, in the first integral in \eq{eq-if5}, we can substitute
$\tilde{\vb}^1-\wb(\tilde{\vb}^0)- \dlt\ub^1 - \tilde{\ub}^1 = \tilde{\vb}^1-\wb(\tilde{\vb}^0) +\wb(\dlt\ub^0) - (\tilde{\ub}^1 + \tilde\wb)$, hence the integral writes
\begin{equation}\label{eq-if10}
  \begin{split}
  \aYs{\tilde{\ub}^1 + \Pibf^{ij}e_{ij}^x(\tilde{\ub}^0)}{\tilde{\vb}^1-\wb(\tilde{\vb}^0) +\wb(\dlt\ub^0) - (\tilde{\ub}^1 + \tilde\wb)}\;. 
\end{split}    
\end{equation}
Now, $\tilde{\vb}^0$ can be chosen such that $\jumpYi{\wb(\dlt\ub^0)-\wb(\tilde{\vb}^0)}{c} = \jumpYi{\tilde\wb}{c}$, therefore the inequality \eq{eq-if10}$\geq 0$ for all $\tilde{\vb}^1$ satisfying \eq{eq-if6a}$_1$ is  equivalent to
\begin{equation}\label{eq-if10a}
  \begin{split}
  \aYs{\tilde{\ub}^1 + \Pibf^{ij}e_{ij}^x(\tilde{\ub}^0)}{\tilde{\vb}^1 - \tilde{\ub}^1}\geq 0\;,\quad \forall \tilde{\vb}^1 \in \Kcal_Y(\nabla\tilde\ub^0)\;,
\end{split}    
\end{equation}
which yields the consistency of \eq{eq-if10a} the assumption on the LVI solution.

Further let us assume:
 $\tilde\ub^0 \in V_0$,  \eq{eq-if10} holds for $(\tilde{\ub}^0,\tilde{\ub}^1)$ and
\begin{equation}\label{eq-if11}
  \begin{split}
 \tilde{\sigmabf}^0(\tilde{\ub}^0,\tilde{\ub}^1):\eebx{\tilde\vb^0} =  \int_\Om \tilde\fb\cdot\tilde\vb^0\;, \quad \forall \tilde\vb^0 \in V_0\;.
\end{split}    
\end{equation}
The last integral on the \lhs in \eq{eq-if5} is identified with the \lhs of the VI in \eq{eq-if12} whereby $ \dlt\tilde{\vb}:= \tilde{\vb}^1 - \tilde{\ub}^1$.  We can show $ \dlt\tilde{\vb}\in \wtilde{\Kcal_Y}(\nabla{\ub}^0) = \Kcal_Y(\nabla{\ub}^0) - \tilde\ub^1$. Indeed, it foll macroscopic contact problems due to \eq{eq-if6}, hence
\begin{equation}\label{eq-if8}
\begin{split}
  \dlt\ub^1 & \in \Kcal_Y(\nabla(\tilde{\ub}^0)) - \wb(\dlt\ub^0) - \tilde{\ub}^1
  = \wtilde{\Kcal_Y}(\nabla{\ub}^0)  - \wb(\dlt\ub^0)\;,\\
  \dlt\tilde{\vb}:= \tilde{\vb}^1 - \tilde{\ub}^1 & \in \Kcal_Y(\nabla(\tilde{\ub}^0))
  - \wb(\dlt\ub^0) - \tilde{\ub}^1 = \wtilde{\Kcal_Y}(\nabla{\ub}^0)  - \wb(\dlt\ub^0) \;.
\end{split}
\end{equation}
Since  $\tilde{\vb}^1 - \tilde{\ub}^1 - \dlt\ub^1 =  \dlt\tilde{\vb} - \dlt\ub^1 \pm  \wb(\dlt\ub^0)$, whereby \eq{eq-if8} yields $\dlt\tilde{\vb} + \wb(\dlt\ub^0) \in \wtilde{\Kcal_Y}(\nabla{\ub}^0)$ and also $\dlt\ub^1 +\wb(\dlt\ub^0) \in \wtilde{\Kcal_Y}(\nabla{\ub}^0)$,
to satisfy \eq{eq-if1}, the increments $(\dlt{\ub}^0,\dlt{\ub}^1)$ must verify \eq{eq-if12}, which proves assertion (i).

(ii)
Assume \eq{eq-if7} holds, \ie the active contact set $\tilde\Gamma_c^*$ is not modified by the increment $(\dlt{\ub}^0,\dlt{\ub}^1)$,
we get
\begin{equation}\label{eq-if7}
\begin{split}
  \jumpYi{\dlt\ub^1 + \nabla\dlt\ub^0\hat y}{*}=0\;,
\end{split}
\end{equation}
thus, $\jumpYi{\dlt\ub^1 + \nabla\dlt\ub^0\hat y}{} = \jumpYi{\tilde\wb}{}= 0$ on $\tilde\Gamma_c^*$, \ie on the true contact boundary:
\begin{equation}\label{eq-if7a}
\begin{split}
 \tilde\Gamma_c^*(x) = \{y \in \Gamma_c|\; g_c^Y(\tilde\ub^1(y),\tilde\ub^0(x,y)) = 0\}.
\end{split}
\end{equation}

As the consequence, \eq{eq-if12}$_1$ becomes a bilateral contact problem. Indeed,
by \eq{eq-if8}$_2$ we have $\jumpYi{\dlt\tilde\vb + \wb(\tilde\vb^0)}{c} = 0$ on $\Gamma_c^*$. Hence, denoting $\dlt\wb := \dlt\tilde{\vb} - \dlt\ub^1$ and recalling $\jumpYi{\tilde{\vb}^0 - \dlt\ub^0}{*} = \jumpYi{\tilde\wb}{*} = 0$, 
due to \eq{eq-if7} we obtain
\begin{equation}\label{eq-if14}
  \begin{split}
\jumpYi{\dlt\wb}{c} = \jumpYi{-\wb(\tilde\vb^0)+\wb(\dlt\ub^0)}{c} = \jumpYi{-\wb(\tilde\vb^0)+\wb(\tilde\vb^0)}{c} = 0\quad \mbox{ on }\Gamma_c^*\;.
  \end{split}    
\end{equation}
Therefore, inequality \eq{eq-if12}$_1$ is equivalent to an equality constrained by the sliding contact set \eq{eq-if16}. 
We find $\dlt\ub^1 \in V_0(\dlt\ub^0,Y_s,x)$ such that
\begin{equation}\label{eq-if15}
  \begin{split}
    \aYs{\dlt\ub^1 + \Pibf^{ij}e_{ij}^x(\dlt{\ub}^0)}{\vb} & = 0\;,\\
 \end{split}    
\end{equation}    
for all $\vb  \in V_0(\bmi{0},Y_s,x)$. The rest of the proof follows by the linearity yielding \eq{eq-if17}-\eq{eq-if20}.
\end{myproof}
  
\begin{myremark}{rem2}
  The limit problem \eq{eq-if12} is nonlinear and couples the pair $(\dlt\ub^1,\dlt\ub^0)$ at the two scales. The linearization arising form (ii) of Proposition~\ref{prop1} is based on the assumption \eq{eq-if7} which was employed in paper \cite{Rohan-Heczko-CaS}. The contact problem is solved locally in ``each microconfiguration'' to satisfy \eq{eq-if10a} for a given macroscopic strain $\eebx{\tilde\ub^0}$. Then the GVE provides the increments  $(\dlt\ub^1,\dlt\ub^0)$ and the updating step follows by \eq{eq-if5a}.
  As explained in the next section, problem \eq{eq-if12} can be solved as macroscopic contact problem while the assumption is weakened to a subset, or dropped fully.
  \end{myremark}
  
  \end{list}

\section{Method of macroscopic contact problem}\label{sec-ma-mcp}

In this section, we present an alternative formulation of the macroscopic problem which arises form the two-scale problem \eq{eq-if12}, Proposition~\ref{prop1} (i), but with a modified, less restrictive assumption on the increments $(\dlt\ub^1,\dlt\ub^0)$.  The two-scale algorithm  proposed in paper \cite{Rohan-Heczko-CaS} is based on the equilibrium equation \eq{eq-if20} which is obtained as the consequence of the fixed bilateral contact, arising from assumption \eq{eq-if7}. 

The newly introduced method of macroscopic contact problem (MC), as explained below, is based on coupled increments $(\dlt\ub^1,\dlt\ub^0)$ due to the bilateral contact on the true contact surfaces $\Gamma_c^*(x)$, however, respects the unilateral contact on $\Gamma_c^\circ(x)\Gamma_c\setminus\Gamma_c^*(x)$ at all micro-configurations $\wtilde\Mcal_Y(x)$, $x\in \Om$.
An iterative algorithm follows the scheme proposed in paper \cite{Rohan-Heczko-CaS}, consisting of alternating local (microscopic) and  global (macroscopic) steps.
The local contact problems \eq{eq-p4} are solved in each microstructure; for given approximation $\tilde\ub^0(x)$, a corrected field $\tilde\ub^1(x,\cdot)$ defined in $Y_s$ is computed which thereby yields the true contact surface $\Gamma_c^*(x) \subset \Gamma_c$. 
The macroscopic (unilateral) contact problems are solved for increments $(\dlt\ub^1,\dlt\ub^0)$ coupled by the linearization based on the assumed bilateral contact on  $\Gamma_c^*(x)$. Updating step follows by virtue of \eq{eq-if5a}.

The MC method is derived from the VI \eq{eq-if1} which governs $(\ub^0,\ub^1) \in W_\#(Y_s,\Om)$. For this, we shall employ some  extra notation which is now introduced.

Let us denote by $(\tilde\ub^0,\tilde\ub^1)$ the reference state characterizing local micro-configurations $\Mcal_Y(x)$ a.e. in $\Om$.
For a given $\wtilde\Mcal_Y(x)$, the contact gap function is given, as
\begin{equation}\label{eq-ma-1}
  \begin{split}
    g_c(\ub^1,\ub^0) & = \jumpY{\Pi_y\ub^0 + \ub^1 - \what y}\;,\\
   \tilde g_c(\dlt\ub^1,\dlt\ub^0) & = \jumpY{\Pi_y(\tilde\ub^0 + \dlt\ub^0) + \tilde\ub^1 + \dlt\ub^1- \what y} = \jumpY{\Pi_y \dlt\ub^0  + \dlt\ub^1 + \tilde\ub^\mic - \what y}\;,
\end{split}
\end{equation}
where $\Pi_y : \ub^0 \mapsto \ol{\ub} = \Pibf^{ij}e_{ij}^x(\ub^0)$ is the affine mapping, generating the homogeneous strain field $\eeby{\ol{\ub}}$ in the whole of $Y$.
For the sake of brevity, we denote by $\tilde\ub^\mic = \tilde\ub^1 + \Pi_y\tilde\ub^0$ the current approximation of the deformation state in $\Mcal_Y(x)$, accordingly we establish the effective stress \ie $\tilde\sigmabf = |Y|^{-1}\int_{Y_s} \Dop \eeby{\tilde\ub^\mic}$.

The two-scale contact set $\Sigma_\Gamma \subset\Om\times\Gamma_c$ reflect the true contact surfaces $\Gamma_c^*(x)$ at $\wtilde\Mcal_Y(x)$,
\begin{equation}\label{eq-ma-2}
  \begin{split}
\Sigma_\Gamma = \{(x,y) \in \RR^d\times\RR^d|\; y \in \Gamma_c\subset\Gamma_c^*(x),\;x\in\Om\}\;.
\end{split}
\end{equation}

The admissibility set for the displacement increments
\begin{equation}\label{eq-mcp3}
\begin{split}
\wtilde\Kcal_\Om^E = \{ \vb \in U_0(\Om) |\;  \wtilde{g}_c^E(\vb) \leq 0 \mbox{ a.e. in } \Sigma_\Gamma\} \;, 
\end{split}
\end{equation}
is defined in terms of the gap function $\what{g}_c^E(\vb)$ which depends locally on the deformed microconfigurations through $\tilde\ub^\mic$ and $\Sigma_\Gamma$, 
see \eq{eq-ma-2},
\begin{equation}\label{eq-mcp4}
\begin{split}
  \wtilde{g}_c^E(\vb(x)) & = \what\Pb:\eebx{\vb} + \wtilde{s}\;,\\
    \what\Pb & :=  \nb\otimes\Delta \what y + \what\Wb_n^E\;,\\
    \wtilde{s} & :=  \nb\otimes\Delta \what y :\eebx{\tilde\ub^0 } + \jumpY{\tilde\ub^1 - \what y}  = \jumpY{\tilde\ub^\mic - \what y}\;,\quad
   \mbox{ a.e. in } \Sigma_\Gamma \;,
\end{split}
\end{equation}
where tensor $\what\Wb_n^E$ is established using the corrector functions computed in \eq{eq-if15}, $(\what W_n^E)_{ij} = \jumpY{\wb^{ij}}$, and $\Delta \what y = \jumpY{\hat y}$. 

We consider $\vb^1 \equiv \ub^1$ and, accordingly, also the perturbations of the solution and the test functions are coupled by \chEM{$\dlt\ub^1 = \dlt\vb^1$}.
Due to the sliding bilateral contact on $\Gamma_c^*$, the micro- and macro-increments are coupled by $\dlt\ub^1 = \what\Wb_n^E\eebx{\dlt\ub^0}$.
With such constraints, \eq{eq-if1} transforms in the following problem: Find $\dlt\ub^0 \in \wtilde\Kcal_\Om^E$ satisfying
\begin{equation}\label{eq-mcp2}
\begin{split}
  \int_\Om\intY_{Y_s} \Dop(\eeby{\tilde\ub^1 + \dlt\ub^1} + \eebx{\ub^0 + \dlt\ub^0}):(\eebx{\vb} - \eebx{\dlt\ub^0})
  \geq & \int_\Om \fb \cdot (\vb - \dlt\ub^0)\;,
\end{split}
\end{equation}
for all $\vb \in \wtilde\Kcal_\Om^E$.

In what follows, by  $\Dop^{E}$ we denote the  
the effective elastic modulus computed according to \eq{eq-if16}-\eq{eq-if19} for locally given actual contact gaps.

To solve \eq{eq-mcp2} numerically, a projection of trial solutions on the admissibility set $\wtilde\Kcal_\Om^E$ is needed. To construct such a projection, we shall consider a saddle point problem involving $(\eb,\ub,\tau)$, with the Lagrangian functional defined, as follows
\begin{equation}\label{eq-mcp5}
\begin{split}
  \Lcal(\eb,\ub,\taubf) & = \Phi(\eb,\ub) + \Psi(\eb,\ub,\taubf)\;,\\
  \Phi(\eb,\ub) & = \frac{1}{2}\int_\Om\Dop^{E}\eb :\eb - \int_\Om \left(\fb\cdot\ub - \tilde\sigmabf :\eb\right)\;,\\
  \Psi(\eb,\ub,\taubf) & = \int_\Om\taubf:(\eb - \eebx{\ub})\;,
\end{split}
\end{equation}
where $\fb:=\tilde \fb + \dlt\fb$ represent the actual load (for the sake of brevity we disregard surface tractions), so that  $\fb\cdot\ub - \tilde\sigmabf :\eb$ is the local out-of balance term, recalling that $\tilde\sigmabf $ is the effective stress defined in the reference configuration.

Clearly, \eq{eq-mcp2} is the necessary condition to be satisfied by $\dlt\ub$ to minimize $\Phi(\eb(\vb),\vb)$ over all admissible increments  $\vb \in \wtilde\Kcal_\Om^E$., However, solutions $\hat\ub:=\dlt\ub$ to problem \eq{eq-mcp2} can be obtained by solving the inf-sup problem for
$(\hat\eb,\hat\ub,\hat\tau)$
satisfying
\begin{equation}\label{eq-mcp6}
\begin{split}
(\hat\eb,\hat\ub,\hat\taubf) = \arg \sup_{\taubf\in\Scal(\Om)} \inf_{\ub \in U(\Om)} \inf_{\eb \in \tilde\Kcal_\Om^E} \Lcal(\eb,\ub,\taubf)\;,
\end{split}
\end{equation}
where $\tilde\Kcal_\Om^E$ is defined in \eq{eq-mcp3} and $\Scal(\Om) = \{\taubf = (\tau_{ij}) \in \SS_2, \tau_{ij} \in L^2(\Om)\}$.

We proceed by the necessary conditions satisfied by solutions of \eq{eq-mcp6}, thus $(\hat\eb,\hat\ub,\hat\taubf) \in \tilde\Kcal_\Om^E \times U(\Om)\times \Scal(\Om)$
must satisfy
\begin{equation}\label{eq-mcp7}
\begin{split}
 (i)& \quad  \int_\Om[\Dop^{E}\hat\eb + \hat\taubf + \tilde\sigmabf]:(\eb - \hat\eb)\geq 0\quad \forall
  \eb \in \tilde\Kcal_\Om^E\;,\\
  (ii)& \quad  \int_\Om \hat\taubf:\eebx{\vb} + \int_\Om \fb\cdot\vb = 0\quad \forall \vb \in U_0(\Om)\;,\\
  (iii)& \quad  \int_\Om  \etabf : (\chEM{\hat\eb} - \eebx{\hat\ub}) = 0\quad \forall\etabf \in \Scal(\Om)\;.
\end{split}
\end{equation}
It is worth to note that the negative multiplier $-\hat\taubf$ expresses the effective macroscopic stress which ensures the equilibrium.
The inequality (i) can be presented as a projection $\hat\eb = \Proj{\hat\db}{\tilde\Kcal_\Om^E}$ of an element  $\hat\db$ of $\SS_2$ on set $\tilde\Kcal_\Om^E$, thus introducing $\hat\db$,
\begin{equation}\label{eq-mcp8}
\begin{split}
& \int_\Om\left(\hat\eb - \hat\db\right):(\eb - \hat\eb)\geq 0\quad \eb \in \tilde\Kcal_\Om^E\;,\\
\mbox{ where }& \hat\db = \hat\eb - \Dop^{E}\hat\eb - \hat\taubf - \tilde\sigmabf\;.
\end{split}
\end{equation}
We shall use a dual representation of the above projection. 
For any $\lam,\vtheta \in L^2(\Sigma_\Gamma)$, let us define
\begin{equation}\label{eq-mcp8a}
\begin{split}
  & \lam_+(x,y) =  \jumpp{\lam}_+ =\chEM{\max\{0,\lam(x,y)\}} \quad \mbox{ a.e. in } {\Sigma_\Gamma}\;,\\
  & \ipprod{\lam}{\vtheta}_{\Sigma_\Gamma} = \int_\Om \int_{\Gamma_c(x)}\lam \vtheta\;,\\
  & \iprod{\etabf}{\eb}_{\Om} = \int_\Om \etabf:\eb\;,
\end{split}
\end{equation}
where above $\Gamma_c(x)$ can be replaced by $\Gamma_c$ if the initial micro-configuration is to be considered, or by $\wtilde\Gamma_c(x)$ for the case of the deformed actual micro-configuration.
Further we introduce the following Lagrangian function $\psi$ and a positive cone $\Ccal_+$,
\begin{equation}\label{eq-mcp9a}
\begin{split}
 \psi(\eb,\lam) & = \frac{1}{2}\nrm{\hat\eb - \hat\db}{\Om}^2
 +  \ipprod{\hat\lam}{\what\Pb:\eb + \wtilde s}_{\Sigma_\Gamma}\;,\\
  \Ccal_+(\Sigma_\Gamma) & = \{\vtheta \in L^2(\Sigma_\Gamma)|\;\vtheta(x,\cdot) \geq 0\;\mbox{ a.e. on }\Gamma_c(x)\}\;.
\end{split}
\end{equation}
Now, problem \eq{eq-mcp8} is equivalent to the saddle point problem involving $\psi$,
\begin{equation}\label{eq-mcp9}
\begin{split}
(\hat\eb,\hat\lam) & = \arg \inf_{\eb \in \Scal(\Om)} \sup_{\lam\in\Ccal_+(\Sigma_\Gamma)} \psi(\eb,\lam)\;,
\end{split}
\end{equation}
The necessary conditions must be satisfied by any solution $(\hat\eb,\hat\lam)$ of \eq{eq-mcp9}, 
\begin{equation}\label{eq-mcp10}
\begin{split}
  \int_\Om(\hat\eb - \hat\db):\eb  + \ipprod{\hat\lambf}{\what\Pb:\eb}_{\Sigma_\Gamma} & = 0\quad \forall \eb \in \Scal(\Om)\;,\\
\chEM{-  \ipprod{\what\Pb:\hat\eb + \wtilde s}{\vtheta-\hat\lam}_{\Sigma_\Gamma} } & \geq 0\quad \forall \vtheta \in \Ccal_+(\Sigma_\Gamma)\;,
\end{split}
\end{equation}
where the inequality imposes a condition which can be expressed by the projection of $\hat\lam - (\what\Pb:\eb + \wtilde s)$ on $\Ccal_+(\Sigma_\Gamma)$
\begin{equation}\label{eq-mcp11}
\begin{split}
  \hat\lam & = \Proj{\chEM{\hat\lam + \what\Pb:\eb + \wtilde s}}{\Ccal_+(\Sigma_\Gamma)} = \jumpp{\chEM{\hat\lam  + (\what\Pb:\eb + \wtilde s)}}_+\;,
\end{split}
\end{equation}
see the definitions in \eq{eq-mcp8a}.

Using the adjoint operator $\what\Pb^*:L^2(\Sigma_\Gamma) \mapsto L^2(\Om)$, such that  $\ipprod{\lam}{\what\Pb:\eb}_{\Sigma_\Gamma} = \iprod{\what\Pb^*\lam}{\eb}_{\Om}$, from \eq{eq-mcp10}$_1$, we can express $\hat\eb$,
\begin{equation}\label{eq-mcp12}
\begin{split}
  \iprod{\Dop^{E}\hat\eb + \hat\taubf + \tilde\sigmabf + \what\Pb^*\hat\lam}{\eb}_{\Om} = 0\quad\forall \eb \in \Scal(\Om)\;,\\
 \mbox{ hence }\quad  \hat\eb = -(\Dop^{E})^{-1}(\hat\taubf + \tilde\sigmabf + \what\Pb^*\hat\lam)\quad\mbox{ a.e. in }\Om\;,
\end{split}
\end{equation}
and substitute in (iii) of \eq{eq-mcp7}, so that the multiplier $\hat\taubf$ can be expressed,
\begin{equation}\label{eq-mcp13}
\begin{split}
  \int_\Om  \etabf : \left( (\Dop^{E})^{-1}(\hat\taubf + \tilde\sigmabf + \what\Pb^*\hat\lam)
  + \eebx{\hat\ub}\right) = 0\quad \forall\etabf \in \Scal(\Om)\;,\\
  \hat\taubf = -\left( \Dop^{E}\eebx{\hat\ub} +\tilde\sigmabf + \what\Pb^*\hat\lam\right)
  \quad\mbox{ a.e. in }\Om\;,
\end{split}
\end{equation}
hence, upon substituting in  \eq{eq-mcp12}$_2$, the strain verifies the consistency with the differentiable displacement field,
\begin{equation}\label{eq-mcp13a}
\begin{split}
 \hat\eb = \eebx{\hat\ub}\quad\mbox{ a.e. in }\Om\;.
\end{split}
\end{equation}
As the next step, $\hat\taubf$ is substituted in (ii) of \eq{eq-mcp7}, so that the macroscopic equilibrium attains the form,
\begin{equation}\label{eq-mcp14}
\begin{split}
  \int_\Om\left(\Dop^{E}\eebx{\hat\ub} + \what\Pb^*\hat\lam\right):\eebx{\vb} & = \int_\Om \fb\cdot\vb -  \int_\Om\tilde\sigmabf:\eebx{\vb}\;, \quad \forall \vb \in U_0(\Om)\;,
\end{split}
\end{equation}
with the \rhs term representing the virtual work of out-of-balance forces (recall $\fb:=\tilde \fb + \dlt\fb$).
Due to \eq{eq-mcp13a} and \eq{eq-mcp8a}, the projection \eq{eq-mcp11} can be rewritten as a nonsmooth equation,
\begin{equation}\label{eq-mcp15}
 \begin{split}
   0 = \chEM{\max\{-\hat\lam,\what\Pb:\eebx{\hat\ub} + \wtilde s\}} \quad \mbox{ a.e. in } {\Sigma_\Gamma}\;.
\end{split}
\end{equation}

\subsection{Solving the macroscopic contact problem}
Equalities \eq{eq-mcp14}-\eq{eq-mcp15} constitute the macroscopic  problem for the increment \chEM{$\dlt\ub^0 = \hat\ub \in  U(\Om)$} involving the multiplier $\hat\lam \in \Ccal_+(\Sigma_\Gamma)$. To find the couple $(\hat\ub,\hat\lam)$, two alternative approaches are straightforward to solve the nonlinear system \eq{eq-mcp14}-\eq{eq-mcp15} iteratively
\begin{list}{}{}
\item \textbf{Tight coupling:} solve \eq{eq-mcp14}-\eq{eq-mcp15} using a linearization of the nonsmooth equation \eq{eq-mcp15}. This approach leads to the use of a nonsmooth version of Newton method after a discretization in space. Note that the discretization of $\Sigma_\Gamma \subset \Om\times\Gamma_c$ leads o high number of DOFs.
\item \textbf{Weak coupling:}  solve \eq{eq-mcp14} and \eq{eq-mcp15} using commuting steps (The Uzawa algorithm): given an approximation of $\lam^k \approx \hat\lam$ at step $k$, solve \eq{eq-mcp14}  for an approximation of $\ub^{k+1} \approx \hat\ub$, then compute a corrected multiplier $\lam^{k+1}$ using a projection step which is based on the following identity arising from \eq{eq-mcp15} with a positive constant $\beta > 0$,
\begin{equation}\label{eq-mcp11a}
\begin{split}
  \hat\lam = \jumpp{\chEM{\hat\lam + \beta (\what\Pb:\eb^{} + \wtilde s)}}_+\;.
\end{split}
\end{equation}

\end{list}

\subsubsection{Application of the Uzawa algorithm}
  The algorithm performs, as follows:

\begin{enumerate}

\item Initiation: for $k = 0$, put $\lam^k = 0$, and set $\beta_0$
  
\item Given  $\lam^k \approx \hat\lam$ at step $k$, solve \eq{eq-mcp14}  to compute  an approximation of $\ub^{k+1} \approx \hat\ub$,

\begin{equation}\label{eq-mcp14a}
\begin{split}
  \int_\Om\Dop^{E}\eebx{\ub^{k+1}}:\eebx{\vb} & = \int_\Om \fb\cdot\vb -  \int_\Om\left(\tilde\sigmabf + \what\Pb^*\lam^k\right):\eebx{\vb}\;, \quad \forall \vb \in U_0(\Om)\;,
\end{split}
\end{equation}

\item Evaluate the local macroscopic strains $\eb^{k+1} = \eebx{\ub^{k+1}}$. (This is an obvious step to emphasize, that the next step can be performed independently point-wise for a. a. $x \in \Om$.)
  
\item Update the multipliers $\lam^k$ using the projection step,

  \begin{equation}\label{eq-mcp11b}
\begin{split}
  \lam^{k+1} = \jumpp{\chEM{\lam^k + \beta_k (\what\Pb:\eb^{k+1} + \wtilde s)}}_+\;,
\end{split}
\end{equation}
  where $\beta_k \in ]0,\bar \beta[$ is a coefficient limited from above by the properties of the elasticity $\Dop^E$ and the operator $\what\Pb$,
 \begin{equation}\label{eq-mcp11c}
\begin{split}
  0< \beta_k < \bar\beta = \min\{2\alpha_\Dop,C_\Pb\}\;,\quad\mbox{ where }\\
  \int_\Om\Dop^{E}\eebx{\vb}:\eebx{\vb} \geq \alpha_\Dop\|\vb\|_{U_0(\Om)}^2\;,\\
  C_\Pb = \max_{(x,y) \in \Sigma_\Gamma, \eb=\eb^T}\{\what\Pb:\eb/|\eb|\}\;.
\end{split}
\end{equation}
Put $k:=k+1$ and go to step 2, unless the convergence $\lam^k \rightarrow \hat\lam$ and $\ub^k  \rightarrow \hat\ub$ is obtained.
 
\end{enumerate}

\subsection{Modifications of the macroscopic contact problem}
The MC method was presented for the two-scale contact surface $\Sigma_\Gamma$ defined in \eq{eq-ma-2} using the complements $\Gamma_c^\circ(x)$ of the actual active sets $\Gamma_c^*(x)$ specific to each $x \in \Om$. The coupling of the increments $\dlt\ub^1$ and $\dlt\ub^0$ by the linearization due to the sliding bilateral contact on $\Gamma_c^*(x)$ presents a kinematic constraint which limits the increment step in the context of the nonlinear unilateral contact problem. The effective stiffness  $\Dop^{H}$  \eq{eq-if16}, being evaluated using the corrector functions $\tilde\wb^{ij}$ computed using \eq{eq-if16}-\eq{eq-if17}, anticipates the active contact on the whole $\Gamma_c^*(x)$ which cannot reduce, but can augment only during the increment.

To alleviate the drawback of the formerly proposed method, we propose its modification which restricts admissible contact set to a neighborhood of $\pd\Gamma_c^*(x)$, such that for a $\gamma$-neighborhood we define
\begin{equation}\label{eq-ma-2*}
  \begin{split}
    \Sigma_\Gamma^\gamma = \{(x,y) \in \RR^d\times\RR^d|\; y \in \Gamma_c^\gamma(x),\;x\in\Om\}\;,\quad
    \Gamma_c^\gamma(x) & = \{y \in \Gamma_c|\; \dist(y, \pd\Gamma_c^*(x)) \leq \gamma \}\;, 
\end{split}
\end{equation}
where  $\gamma>0$ is small.
Accordingly, the  stiffness  $\Dop^{H}$ employed in \eq{eq-mcp14} should be computed using the bilateral sliding contact on the reduced active set
$\hat\Gamma_c^*(x) := \Gamma_c^*(x)\setminus \Gamma_c^\gamma(x)$. This influences the corrector fields $\hat\Wb^{ij}$ and consequently the operator $\what\Pb$ defined using \eq{eq-mcp4}, correspondingly. The modified  macroscopic contact problem \eq{eq-mcp14}-\eq{eq-mcp15} involves the restricted contact set $\Sigma_\Gamma^\gamma$. Its discretization leads to a problem with a low number of DOFs associated with the $\gamma$-neighborhood of $\pd\Gamma_c^*(x)$, when compared to the former method requiring the discretized set $\Sigma_\Gamma$. This enables for an efficient solving of the  modified system  \eq{eq-mcp14}-\eq{eq-mcp15} using the ``tight coupling'' alternatively, rather then using the Uzawa algorithm.

\section{Alternative formulations of the microscopic problem}\label{sec-micro}
The aim of this section is to introduce a dual formulation of the LVI \eq{eq-p4} interpreted for the increments \eq{eq-if5a}. Let us recall the limit variational inequality \eq{eq-if1} and the definition of the admissibility sets $W_\#$ in \eq{eq-p3b}.  The incremental field $\dlt\ub^1 = \ub^1 - \tilde\ub^1$, where $\ub^1$ is the unknown microscopic displacement, while considering fixed $\ub^0\approx\tilde\ub^0$, must be consistent with the limit VI \eq{eq-if1}. Upon substituting there for $\ub^1$, we get the following problem for $\dlt\ub^1 \in W_\#(Y_s,\Om)$ satisfying
\begin{equation}\label{eq-mic0}
\begin{split}
  \int_\Om\intY_{Y_s} \Dop(\eeby{\dlt\ub^1 + \tilde\ub^1} + \eebx{\tilde\ub^0}):(\eeby{\vb^1} + \eebx{\vb^0} - \eeby{\dlt\ub^1 + \tilde\ub^1} -\eebx{\tilde\ub^0})
  \geq & \int_\Om \bar \fb \cdot (\vb^0 - \tilde\ub^0)\;,
\end{split}
\end{equation}
for all $(\vb^0,\vb^1) \in W_\#(Y_s,\Om)$. 
Further some new notation will be used for the sake of brevity. Let us define $\tilde\sb  := \what y - \tilde\ub^\mic$ and 
\begin{equation}\label{eq-mic0a}
  \begin{split}
    \wtilde\Kset_Y & = \{ \vb \in \Hpdbav(Y_s)\;|\; \wtilde g_c(\vb)\leq 0 \mbox{ a.e. on }\Gamma_c\} = \Kcal_Y(\nabla\tilde\ub^0) - \tilde\ub^1\;,\\
    \mbox{ where } \quad
  \wtilde g_c(\vb) & = \jumpY{\vb + \tilde\ub^\mic - \what y} = \jumpY{\vb - \tilde\sb}
  =\wtilde g_c'(\vb  - \tilde\sb ) = \wtilde g_c'(\vb) - \wtilde g_c'(\tilde\sb)\;,\\ 
\end{split}
\end{equation}
so that $\wtilde g_c'(\vb) = \jumpY{\vb}$ is linear. 
We shall apply the following substitutions which are coherent with fixing the macroscopic test functions $\vb^0$ without any loss of generality,
\begin{equation}\label{eq-mic2}
\begin{split}
  \vb^1 & := \tilde\ub^1 + \tilde\vb\;,\\
  \vb^0 & := \tilde\ub^0\;.
\end{split}
\end{equation}
Now \eq{eq-mic0} reads: find $\dlt\ub^1 \in \wtilde\Kset_Y$, such that
\begin{equation}\label{eq-mic1}
\begin{split}
  \intY_{Y_s} \Dop(\eeby{\dlt\ub^1 + \tilde\ub^1} + \eebx{\tilde\ub^0}):(\eeby{\tilde\vb}  - \eeby{\dlt\ub^1})
  \geq & 0\;,\quad \mbox{ a.e. in }\Om,
\end{split}
\end{equation}
for all $\tilde\vb \in  \wtilde\Kset_Y$.

We proceed by introducing the Lagrangian function associated with the above variational inequality, and some further notation,
\begin{equation}\label{eq-mic3}
\begin{split}
  \Lcal_Y(\ub,\lam) & = \frac{1}{2}\aYs{\ub}{\ub} + \intY_{Y_s}\tilde\sigmabf^\mic:\eeby{\ub} + \intY_{\Gamma_c}\lam\wtilde g_c(\ub)\;,\\
  \Ccone_+(\Gamma_c) & = \{\vtheta \in L^2(\Gamma_c)\;|\;\vtheta \geq 0 \mbox{ a.e. on }\Gamma_c\}\;,\\
  \tilde\sigmabf^\mic & = \Dop\eeby{\tilde\ub^\mic}\;,
\end{split}
\end{equation}
where $\Ccone_+(\Gamma_c)$ is the convex cone, $\lam$ denotes the Lagrange multiplier, and $\tilde\sigmabf^\mic$ is the actual stress in the particular micro-configuration $\Mcal_Y(x)$ located at $x \in \Om$; Note that $\tilde\sigmabf^\mic$ is the only data relating the actual macroscopic response to $\Mcal_Y(x)$ which is characterized by the fluctuating field $\tilde\ub^1(x,\cdot) \in \Hpdbav(Y_s)$.
We consider the following saddle point problem:
\begin{equation}\label{eq-mic4}
\begin{split}
(\hat\ub,\hat\lam) = \arg \sup_{\lam \in \Ccone_+(\Gamma_c)} \inf_{\ub \in \Hpdb(Y_s)}  \Lcal_Y(\ub,\lam) \;.
\end{split}
\end{equation}
Solutions $(\hat\ub,\hat\lam)$ must satisfy the following necessary conditions,
\begin{equation}\label{eq-mic5}
\begin{split}
  \aYs{\hat\ub}{\vb} + \intY_{\Gamma_c}\hat\lam\wtilde g_c'(\vb) & = -\intY_{Y_s}\tilde\sigmabf^\mic:\eeby{\vb} \;,\quad \forall \vb \in \Hpdb(Y_s)\;,\\
  \intY_{\Gamma_c}\wtilde g_c(\hat\ub) (\vtheta - \hat\lam) & \chER{\leq} 0\;,\quad \forall \vtheta \in \Ccone_+(\Gamma_c)\;.
\end{split}
\end{equation}
The second conditions yields the projection relationship for $\hat\lam$,
\begin{equation}\label{eq-mic6}
\begin{split}
  \hat\lam = \Proj{\hat\lam \chER{+} \beta\wtilde g_c(\hat\ub)}{\Ccone_+(\Gamma_c)}\;,\\
  \Rightarrow 
  \max\{-\lam,\, \beta \wtilde g_c(\hat\ub)\} = 0\;,\quad \mbox{ a.e. on }\Gamma_c\;,
\end{split}
\end{equation}
where $\beta > 0$ is arbitrary real strictly positive parameter. Projection \eq{eq-mic6}
can be employed in an implementation of the Uzawa algorithm to solve \eq{eq-mic5}. In this case, $\beta$ is bounded from above according constants arising from the ellipticity of $\aYs{\cdot}{\cdot}$ and the Lipschitz continuity of the gap function $\wtilde g_c$ characterizing the constraint functional.

Another way of solving  \eq{eq-mic5} is to write \eq{eq-mic6} as a nonsmooth equality and to express $\ub$ using the resolvent operator $\Acal^{-1}$ associated with the bilinear form  $\aYs{\cdot}{\cdot}$. For this, let us denote by $\Acal,\Gcal,$ and $\Bcal$ the linear operators, such that
\begin{equation}\label{eq-mic7}
\begin{split}
  \Acal:\Hpdbav(Y_s) & \mapsto \Hpdbav(Y_s)^*\;,\quad \aYs{\hat\ub}{\vb} = \ipYs{\Acal\ub}{\vb}\;,\\
  \Bcal:L^2(Y_s;\SS_2) & \mapsto \Hpdbav(Y_s)^*\;,\quad \intY_{Y_s}\sigmabf:\eeby{\vb} = \ipYs{\Bcal\sigmabf}{\vb}\;,\\
\Gcal:\Hpdbav(Y_s) & \mapsto :L^2(\Gamma_c)\;,\quad \intY_{\Gamma_c}\wtilde g_c'(\ub) \vtheta
= \ipGc{\Gcal\ub}{\vtheta} = \ipYs{\Gcal^*\vtheta}{\ub}\;,
\end{split}
\end{equation}
where $\Gcal^*$ is the adjoint operator to $\Gcal$, space $L^2(Y_s;\SS_2) = \{\eb = (e_{ij}),\;e_{ij} \in L^2(Y_s),\; \eb \in \SS_2\}$ and $\Hpdbav(Y_s)^*$ represent the dual space of $\Hpdbav(Y_s)$. Note that $\wtilde g_c(\ub) = \wtilde g_c'(\ub - \wtilde\sb) = \Gcal(\ub- \wtilde\sb)$.

With this notation in hand, problem \eq{eq-mic5} can be rewritten, as follows (we drop the hat $\hat{}$, thus $(\ub,\lam)$ denotes the solutions)
\begin{equation}\label{eq-mic8}
\begin{split}
  \ipYs{\Acal\ub}{\vb} + \ipYs{\Bcal\tilde\sigmabf^\mic}{\vb} + \ipYs{\Gcal^*\lam}{\vb} & = 0\;,\quad \forall \vb \in \Hpdbav(Y_s)\;,\\
  \min\{\chER{-}\ipGc{\Gcal(\ub - \tilde\sb)}{\vtheta}|\;\ipGc{\lam}{\vtheta}\} & = 0\;, \quad \forall \vtheta \in \Ccone_+^*(\Gamma_c)\;,
\end{split}
\end{equation}
where $\Ccone_+^*(\Gamma_c) = \Ccone_+(\Gamma_c)$ denotes the self-dual cone of $\Ccone_+(\Gamma_c)$ (so the notation introduced for formal reasons). Since $\tilde\sigmabf^\mic$ can be expressed in terms of $\tilde\ub^\mic$, see \eq{eq-mic0a},
\begin{equation}\label{eq-mic8a}
\begin{split}
\ipYs{\Bcal\tilde\sigmabf^\mic}{\vb} =  \ipYs{\bar\Acal\tilde\ub^\mic}{\vb}\;,\quad  \bar\Acal:\Hdb(Y_s) & \mapsto \Hpdbav(Y_s)^*\;,
\end{split}
\end{equation}
however, to keep the treatment general enough, we keep using operator $\Bcal$.
Due to the the strong monotonicity a Lipschitz continuity of $\Acal$, there exists an inverse operator $\Acal^{-1}:\Hpdbav(Y_s)^*\mapsto\Hpdb(Y_s)$, such that $\ub = -\Acal^{-1}(\Gcal^*\lam + \Bcal\tilde\sigmabf^\mic)$. Upon substitution in \eq{eq-mic8}$_2$, this identity yields
\begin{equation}\label{eq-mic9}
\begin{split}
   \max\{\chER{-}\ipGc{\Gcal(\tilde\sb + \Acal^{-1}(\Gcal^*\lam + \Bcal\tilde\sigmabf^\mic))}{\vtheta}|\,-\ipGc{\lam}{\vtheta}\} & = 0\;, \quad \forall \vtheta \in \Ccone_+^*(\Gamma_c)\;,
\end{split}
\end{equation}
which can be written point-wise, due to the property $\Ccone_+^*(\Gamma_c) = \Ccone_+(\Gamma_c)$,
\begin{equation}\label{eq-mic10}
\begin{split}
  \chER{\min\{ \Ccal\lam + h|\,\lam\}} & = 0\;, \mbox{ a. e. on } \Gamma_c\;,\\
    \Ccal & = \Gcal\Acal^{-1}\Gcal^*\;,\\
    h & = \Gcal\Acal^{-1}\Bcal\tilde\sigmabf^\mic + \Gcal\tilde\sb\;,
\end{split}
\end{equation}
With reference to \eq{eq-mic8a}, note that $\tilde\ub^\mic$ is not $Y$-periodic, in general,  being given by the projection, since $\Acal^{-1}\bar\Acal$ is not an  
identity operator.

Below we introduce an approximation of$\Gcal$ using a symmetric definition of the contact gap function. By virtue of the numerical discretization of \eq{eq-mic8} using the finite element method, formulation \eq{eq-mic10} relies on the inversion a stiffness matrix $\Abm$ corresponding to operator $\Acal$. Assuming a periodic microstructure, computing of $\Abm^{-1}$ is rather efficient in the context of the whole two-scale algorithm, since $\Abm^{-1}$ is shared by all microconfigurations $\Mcal_Y(x)$.

\subsection{A symmetric approximation of the contact conditions}

For numerical solving problem \eq{eq-mic10}, an approximation of the contact  function is need. For this, operator $\Gcal$ can replaced by its approximation $\wtilde\Gcal$ which takes into account evaluation of the contact gap using a symmetric linearization of the contact constraint. For this we need to introduce homologous pairs of points on the matching boundaries $\Gamma_+ \subset \Gamma_c$ and $\Gamma_- \subset \Gamma_c$, recalling $\Gamma_+ \cup \Gamma_- = \Gamma_c$.

Let $y^+ \in \Gamma_+$, and define its contact-homologous point $\tilde y^- = \psibf_+^-(y^+) := \xi\nb(y^+) + y^+$, such that $\tilde y^- \in \Gamma_-$ for some $\xi \in \RR$. Reciprocally, we introduce $\tilde y^+ = \psibf_-^+(y^-) := \xi\nb(y^-) + y^-$, such that $\tilde y^+ \in \Gamma_+$ for some $\xi \in \RR$. We call $y^\pm$ the ``master'' points, whereas $\tilde y^\mp$ are called the ``slave'' points. These pairs $(y^+,\tilde y^-)$ and $(y^-,\tilde y^+)$ of homologous points enable us to define
a symmetric gap function $\wtilde g_\pm(\ub) \approx \wtilde g_c'(\ub)$. For this, we yet need homologous master points: we consider a scalar parameter $t \in T_c = [0,1]^{d-1}$, $d = 2,3$ (for 2D and 3D problems, respectively) and two bijective mappings $\gamma_+:t\mapsto y^+ \in \Gamma_+$ and $\gamma_-:t\mapsto y^- \in \Gamma_-$. In this way, the homologous master pairs $\gamma_\pm(t) = (\gamma_+(t),\gamma_-(t)) = (y^+,y^-)$ are established for any $t \in T_c$, and employed to define the symmetric approximation of the gap function,
\begin{equation}\label{eq-mic11}
\begin{split}
  \wtilde g_\pm(\ub,t) & := \frac{1}{2}(g_+(\ub) + g_-(\ub))\;,\quad \mbox{ for a.e. } t \in T_c\;,\\
  \mbox{ with }\quad g_+(\ub) & = \nb(y^+)\cdot\left(\ub(y^+) - \gammabf^-(\ub,\tilde y^-)\right)\;,\\
  g_-(\ub) & = \nb(y^-)\cdot\left(\ub(y^-) - \gammabf^+(\ub,\tilde y^+)\right)\;,
\end{split}
\end{equation}
where $\gammabf^+(\ub,\tilde y^+)$ is a linear mapping which serves an approximation of displacement $\ub(\tilde y^+)$ at $\tilde y^+$, the contact-homologous point associated to the master point $y^-$ by mapping $\psibf_-^+(y^-)$. In analogy, $\gammabf^-(\ub,\tilde y^-)$ is introduced.
It is wort to emphasize, that  $\wtilde g_\pm(\ub,t)$ is defined ``point-wise'' \wrt the contact surface parametrization $t$.


The above arrangements allow us to define the following approximation $\wtilde\Gcal\approx \Gcal$,
\begin{equation}\label{eq-mic12}
\begin{split}
  \wtilde\Gcal\ub = \wtilde g_\pm(\ub,t)\;, \mbox{ for a.e. } t \in T_c\;.
\end{split}
\end{equation}
While the construction of $\wtilde\Gcal$ is straightforward, being given by the chain of relationships \eq{eq-mic11}, its adjoint $\wtilde\Gcal^*$ is defined rather in the implicit way, though given by the same relationships by virtue of the collocation at $t \in T_c$. For the  contact conditions  approximated by virtue of \eq{eq-mic11}-\eq{eq-mic12}, the formulation \eq{eq-mic10} is modified straightforwardly.


\section{Numerical examples}\label{sec-numex}

\def\Msquare{\mathcal{M}_\square}
\def\Mcircle{\mathcal{M}_\circ}

To illustrate the numerical solutions of the nonlinear two-scale problem treated in this paper, we report three numerical examples tests. 
First we test the microroblem solutions on two microstructure geometries, then we consider the global, two-scale problem and show the performance of the two proposed methods used to comlute response to the uniaxial compression loading,  and finally we report a simulation of a heterogeneous, fissured short cantilever subject to bending.

Material properties of the solid material are given by the Young modulus $E=2.3\unit{GPa}$ and the Poisson ratio $\nu=0.3$. We confine to 2D problems under the assumption of plane strain.

\subsection{Local problems}

\begin{figure}
  \centering
  \begin{minipage}{.48\textwidth}
    \includegraphics[width=\textwidth]{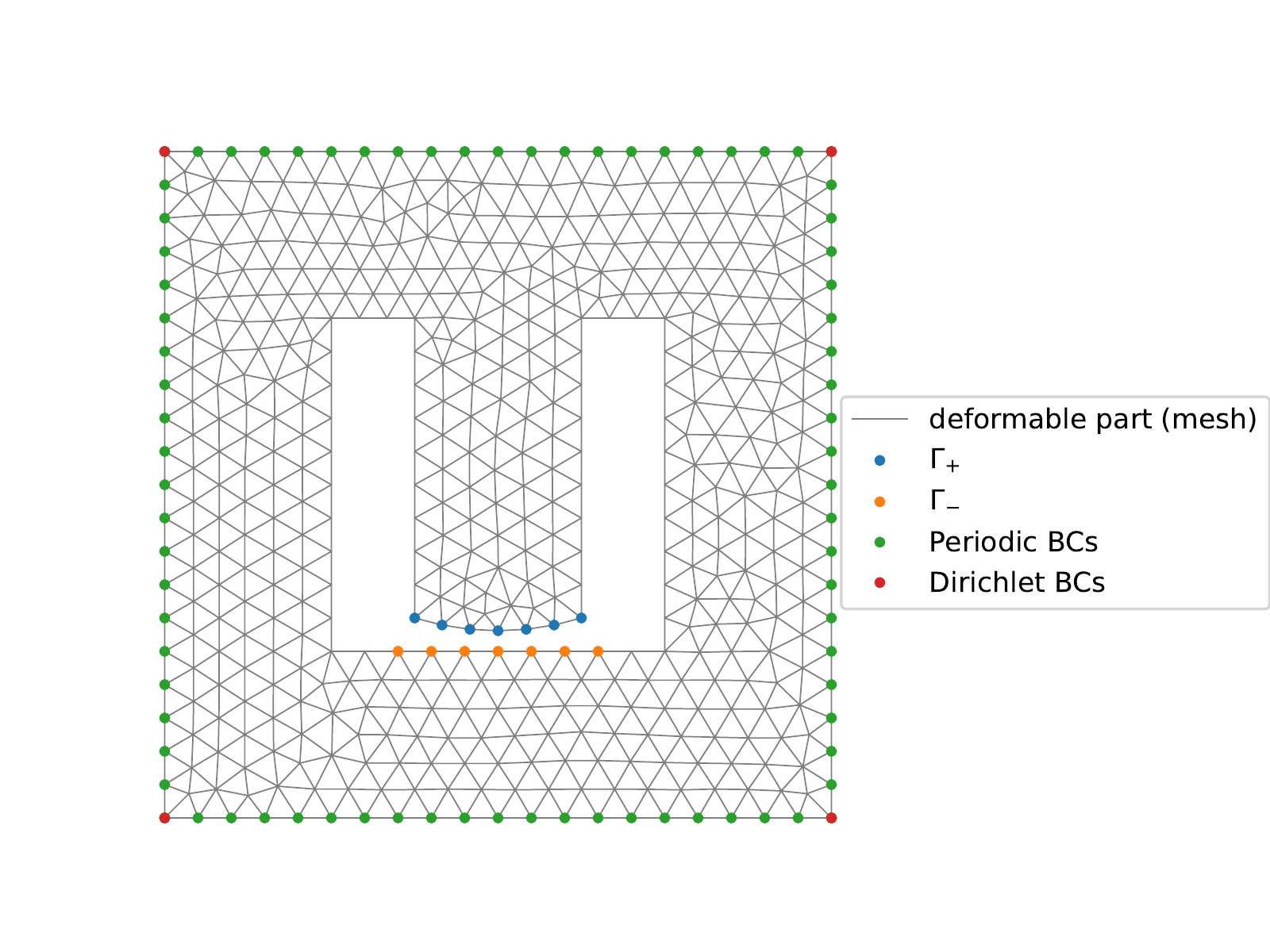}
    \caption{FEM mesh of the ,,uniaxial'' microscopic periodic cell ($\Msquare$).}
    \label{fig-micro-cell-uniaxial}
  \end{minipage}
  \hfill
  \begin{minipage}{.48\textwidth}
    \centering
    \includegraphics[width=\textwidth]{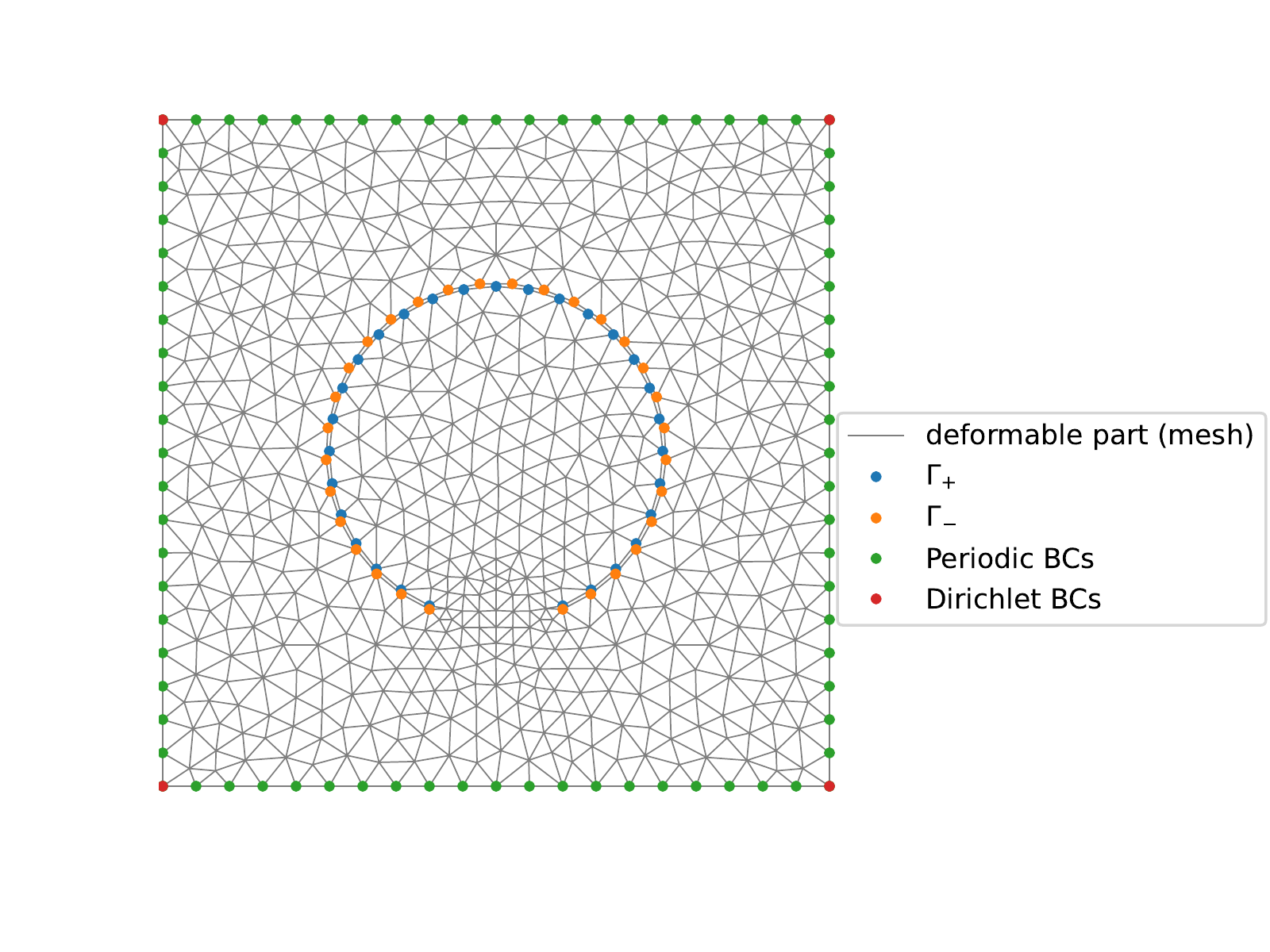}
    \caption{FEM mesh of the microscopic periodic cell with circular inclusion ($\Mcircle$).}
    \label{fig-micro-cell-circle}
  \end{minipage}
\end{figure}

The solution of local problems is illustrated using three examples.
Two different microstructures were considered -- $\Msquare$, shown in Fig.~\ref{fig-micro-cell-uniaxial}, and $\Mcircle$, Fig.~\ref{fig-micro-cell-circle} -- and two different modes of deformation were prescribed:
\begin{equation}
  \mathbf{e}^0_\mathrm{A} =  0.04 \, \begin{bmatrix}
    0.35 & 0 \\
    0 & -1 \\
  \end{bmatrix} \; , \quad
  \mathbf{e}^0_\mathrm{B} =  0.05 \, \begin{bmatrix}
    0 & 1 \\
    1 & 0 \\
  \end{bmatrix} \; .
\end{equation}
$\Msquare$ was subjected to $\mathbf{e}^0_\mathrm{A}$ and $\Mcircle$ was subjected to both $\mathbf{e}^0_\mathrm{A}$ and $\mathbf{e}^0_\mathrm{B}$.


\begin{figure}
  \centering
  \includegraphics[width=.32\textwidth]{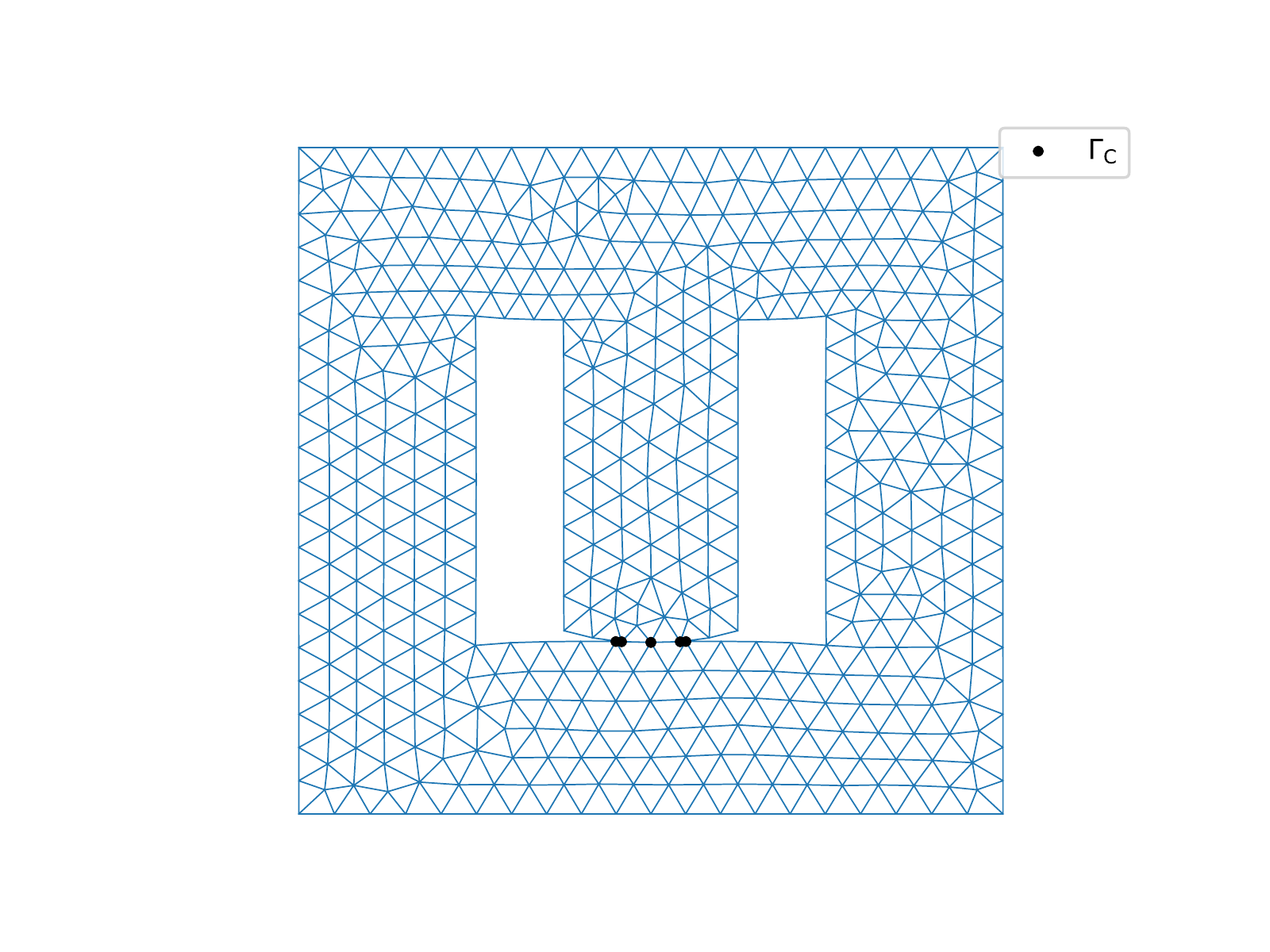}
  \includegraphics[width=.32\textwidth]{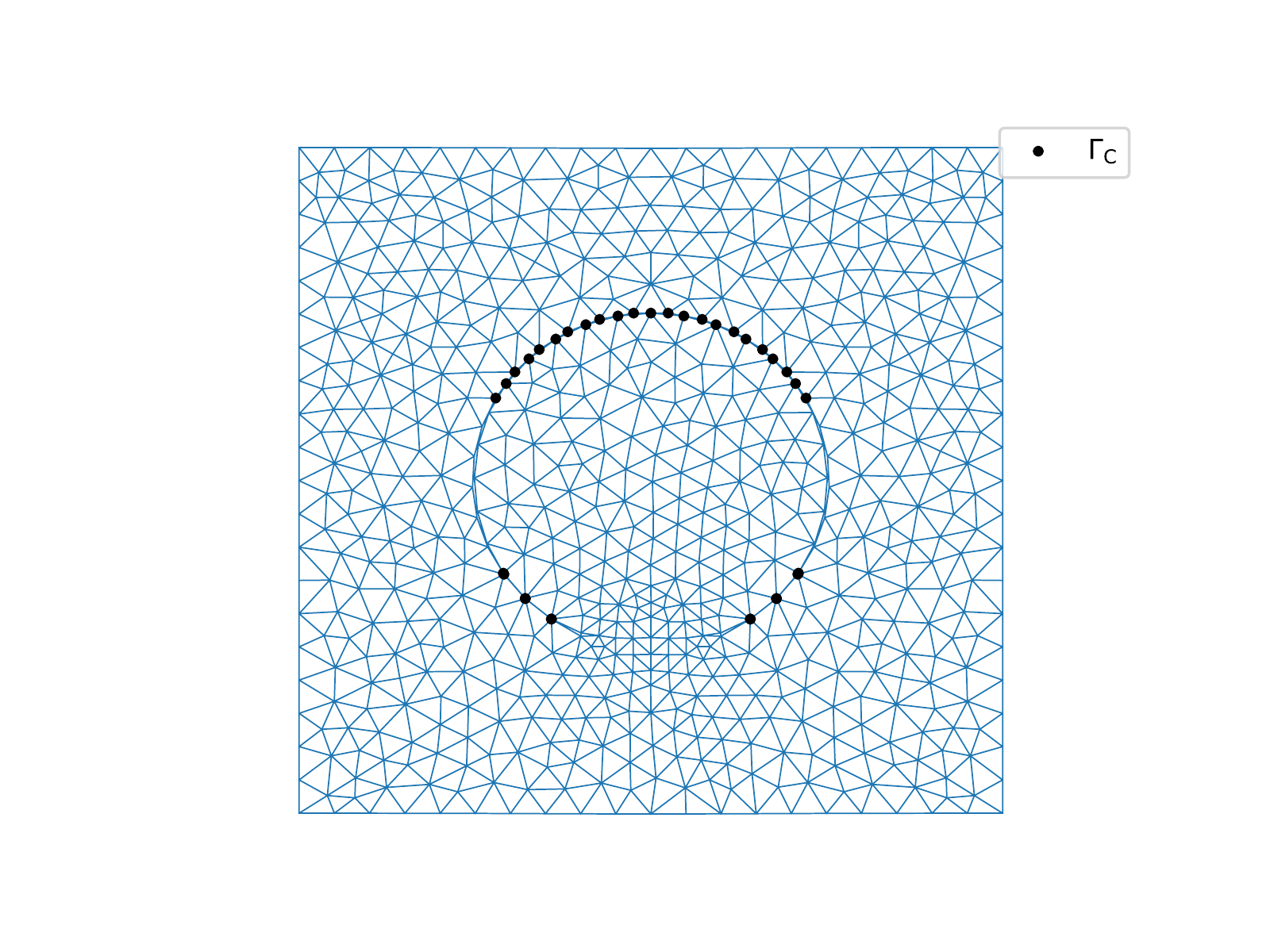}
  \includegraphics[width=.32\textwidth]{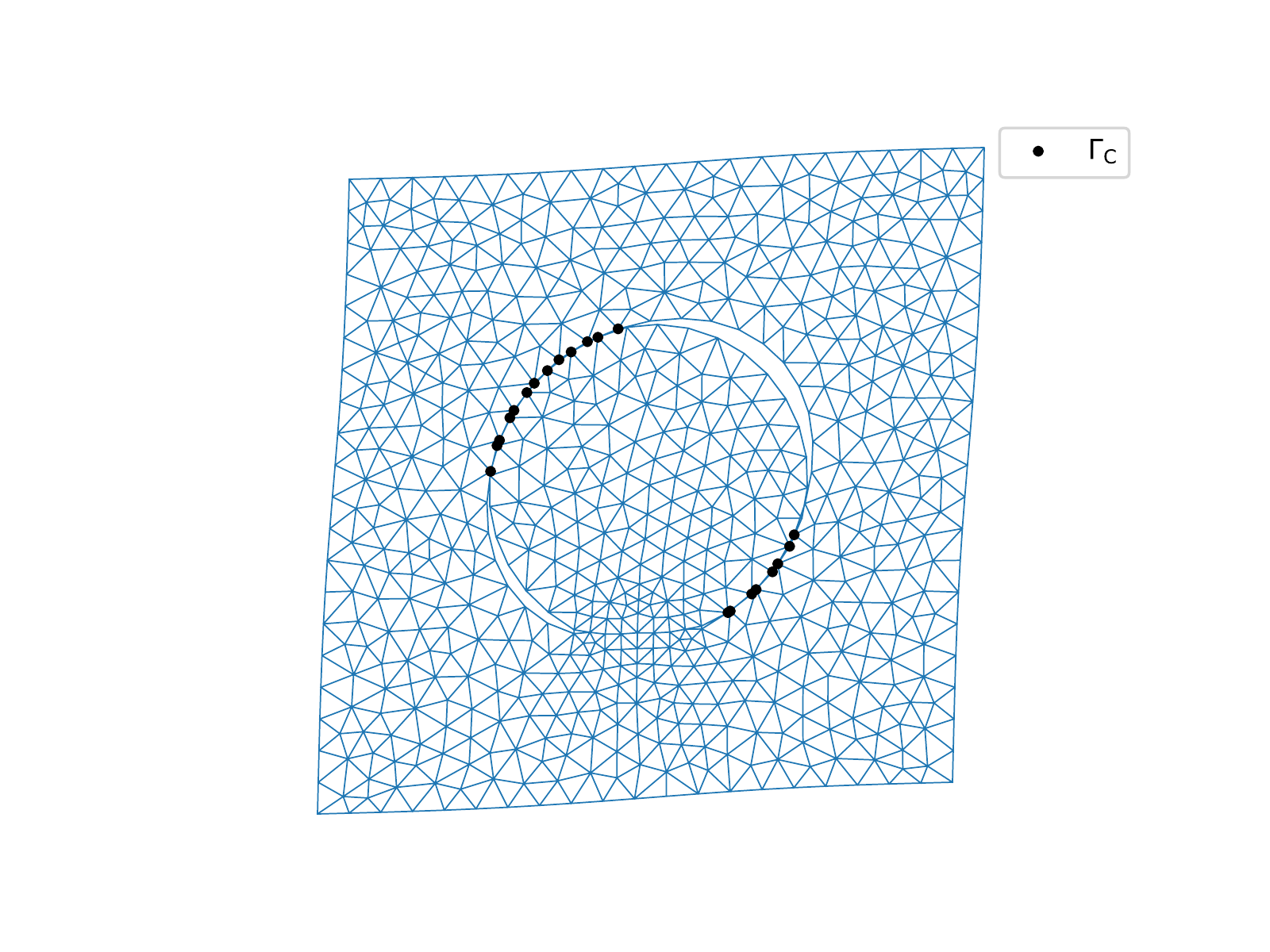}
  \includegraphics[width=.32\textwidth]{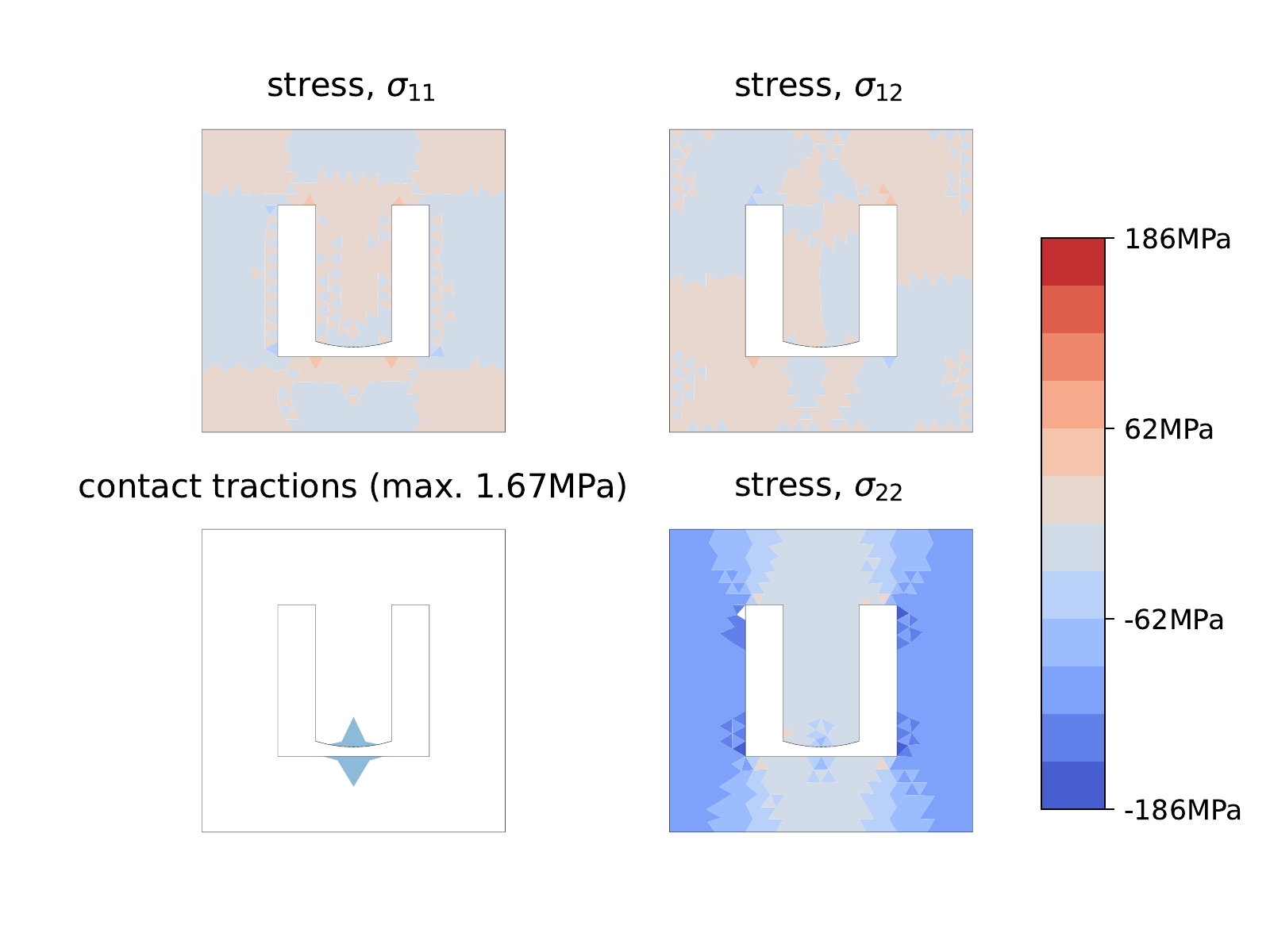}
  \includegraphics[width=.32\textwidth]{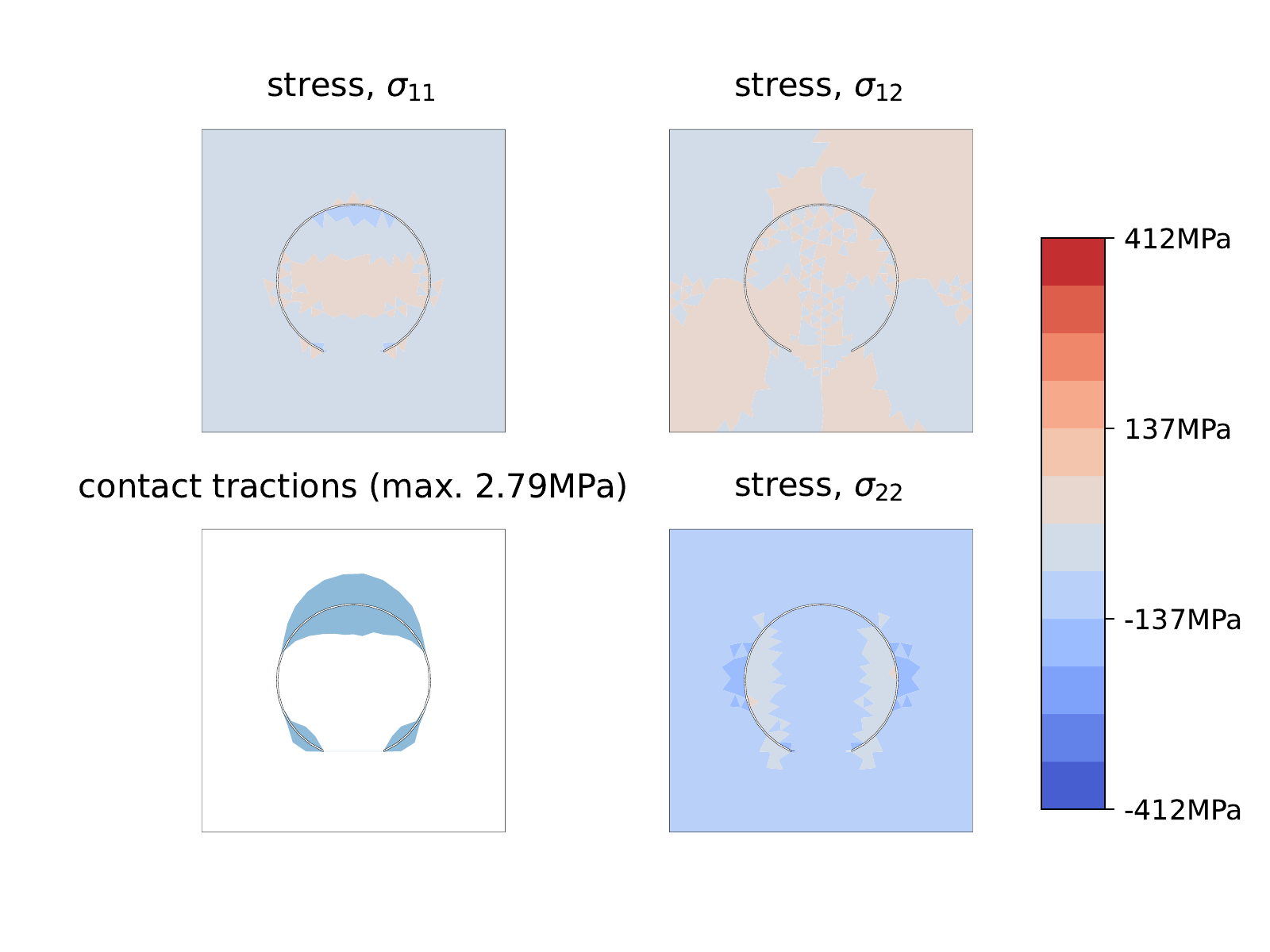}
  \includegraphics[width=.32\textwidth]{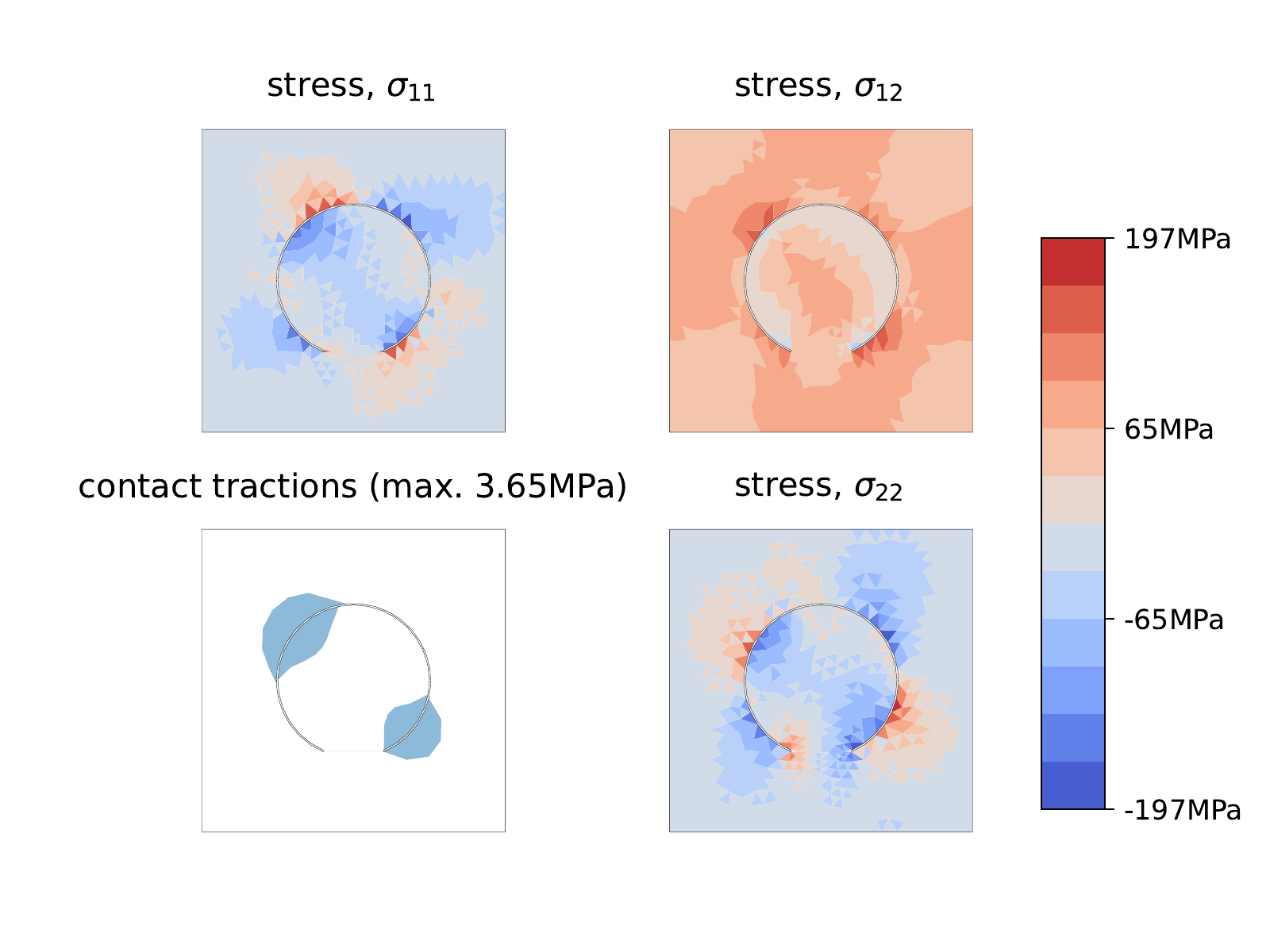}
  \caption{
    Solutions of microscopic (local) problems.
    Deformed shapes of the microscopic cell (top row) and stress components (bottom row).
    Left: $\Msquare$, $\mathbf{e}^0_\mathrm{A}$,
    middle: $\Mcircle$, $\mathbf{e}^0_\mathrm{A}$,
    right: $\Mcircle$, $\mathbf{e}^0_\mathrm{B}$.
  }
  \label{fig-micro-solutions}
\end{figure}

\begin{figure}
  \centering
  \begin{minipage}{.48\textwidth}
    \centering
    \includegraphics[width=\textwidth]{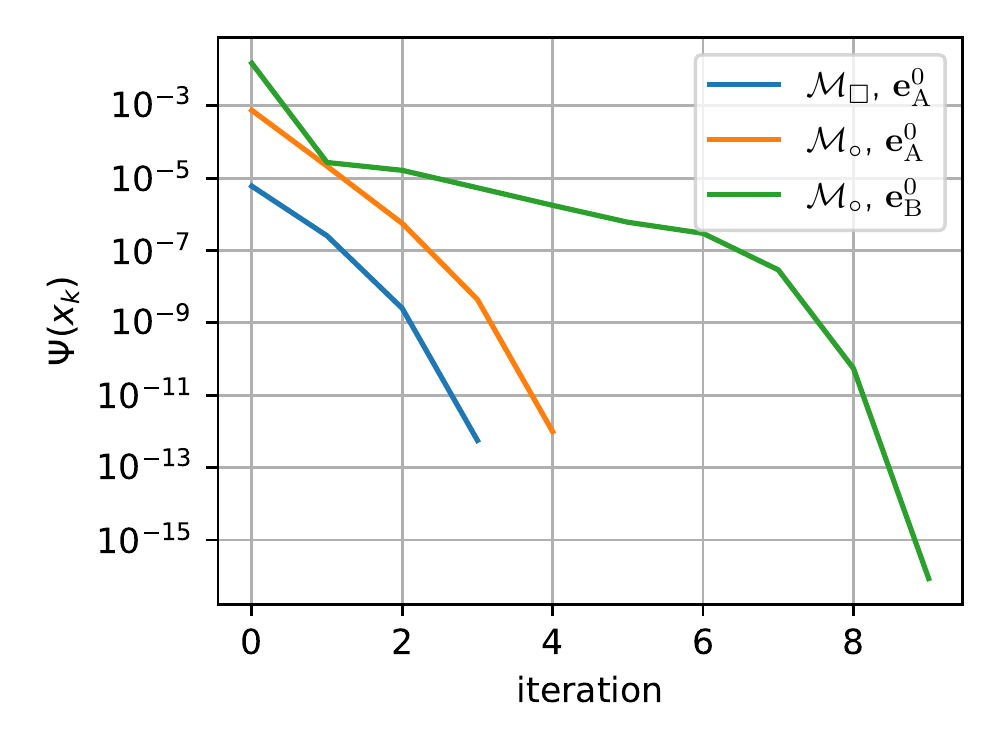}
    \caption{Convergence of local problems.}
    \label{fig-micro-convergence}
  \end{minipage}
  \hfill
  \begin{minipage}{.48\textwidth}
    \centering
    \includegraphics[width=\textwidth]{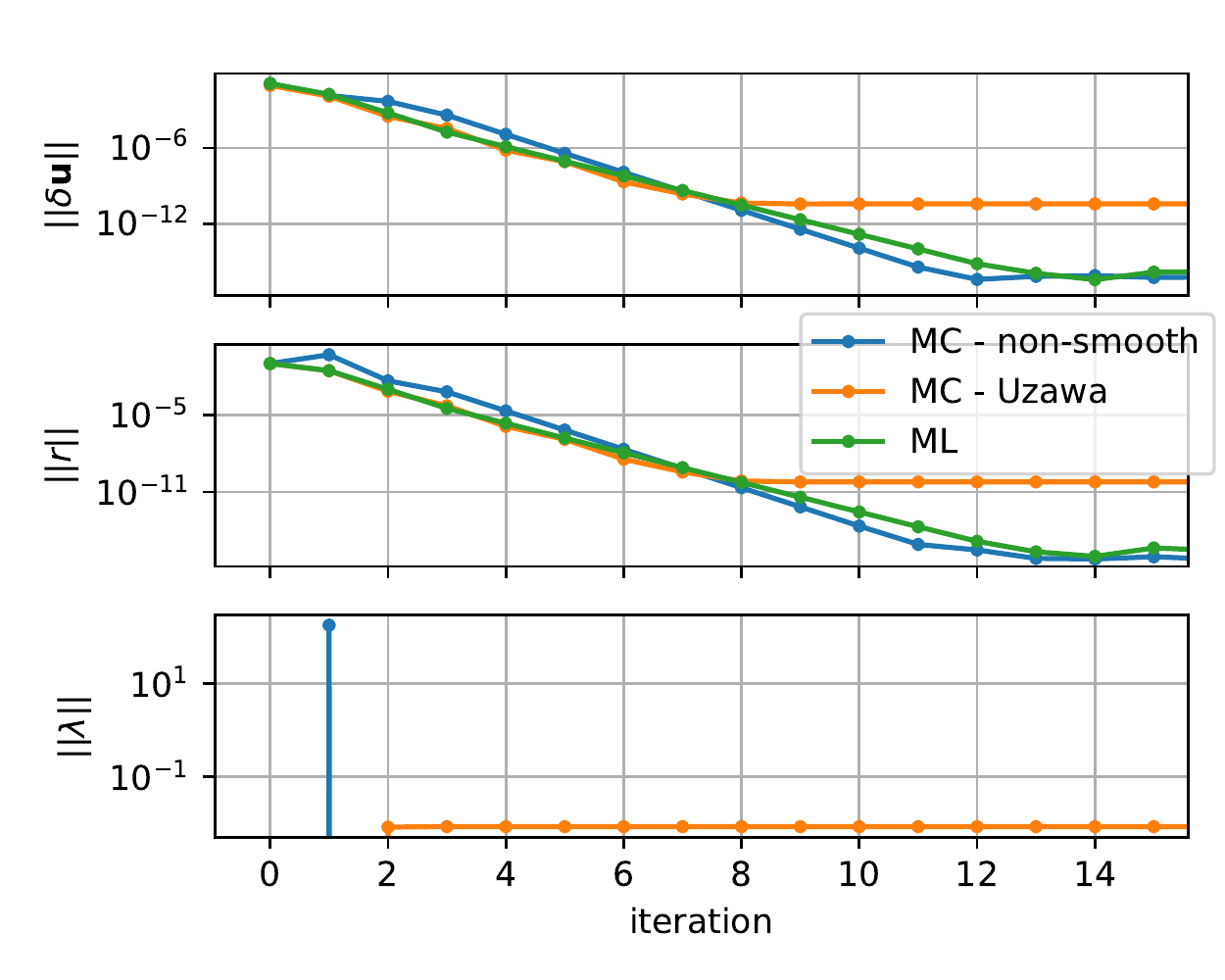}
    \caption{Convergence of the global algorithms in the case of uniaxial compression. Multipliers $\lam$ overlap for the two versions of the MC method.}
    \label{fig-macro-convergence-all}
  \end{minipage}
\end{figure}

In Fig.~\ref{fig-micro-solutions}, shows for each of the above introduced local problems:
the deformed mesh with true contact boundary $\Gamma_c^*$, stress components, and contact tractions are depoicted for the three local problems solved. 
Fig.~\ref{fig-micro-convergence} shows convergence of the non-smooth solver based on the semismooth Newton method, see \cite{Rohan-Heczko-CaS} and \cite{DeLuca-Facchinei-Kanzow-MathProg1996}.

\subsection{Uniaxial compression -- comparison of global algorithms}

In this example, we compare three methods to solve the global problem: The macroscopic contact (MC) method implemented in two variants based on a) the non-smooth solver, or b) on the Uzawa algorithm, is compared  the macroscopic linear (ML) method.
The macroscopic domain is a unit square, $x \in \Om = [0, 1] \times [0, 1]$, meshed by two four-node quadrilateral elements with bilinear approximation of displacement.
The following conditions enforce homogeneous strain distribution:
\begin{equation}
  u_1 = 0 \quad \text{on } \Gamma_\mathrm{left} \; , \quad
  u_2 = 0 \quad \text{on } \Gamma_\mathrm{bottom} \; , \quad
  u_1 = \mathrm{const.} \quad \text{on } \Gamma_\mathrm{right} \; , \quad
  u_2 = \mathrm{const.} \quad \text{on } \Gamma_\mathrm{top} \; .
\end{equation}
Uniform pressure $\bar t = 0.1\unit{GPa}$ was prescribed on $\Gamma_\mathrm{top}$.
The microscopic domain $Y$ is shown in Fig.~\ref{fig-micro-cell-uniaxial}.

Fig.~\ref{fig-macro-convergence-all} shows the norms of selected quantities (the displacement increments -- corrections, $\dlt\ub$, the out-of balance residual, $\rb$, and the $\lam$ multiplier) related to the iterative solution of the global problem.
All three methods eventually hit a convergence plateau; the Uzawa algorithm stops its progress at $\| \rb\|$ around $10^{-10}$, while the other two methods reach values of $10^{-16}$.
  It can be seen that the ``MC -- non-smooth'' variant starts slightly worse, but after a few iterations achieves faster rate of convergence.
  Also, the multipliers $\hat\lam$ in case of the ``MC -- non-smooth'' variant are non-zero only at iteration 1 due to the varying nature of $\Gamma_c^\circ$. 
  
\subsection{Cantilever bending}

The macroscopic domain is a unit square, $x \in \Om = [0, 1] \times [0, 1]$, meshed by an array of $4\times 4$ four-node quadrilateral elements with bilinear approximation of displacement.
All displacements are fixed at the bottom edge, $\Gamma_\mathrm{bottom}$, and a horizontal traction $\bar t = 0.01\unit{GPa}$ is prescribed at the top edge, $\Gamma_\mathrm{top}$.
$\Mcircle$ was used as the microscopic cell (see Fig.~\ref{fig-micro-cell-circle}) and the material properties were the same as in the previous examples.

\begin{figure}
  \centering
  \begin{tabular}{cc}
  \includegraphics[width=.49\textwidth]{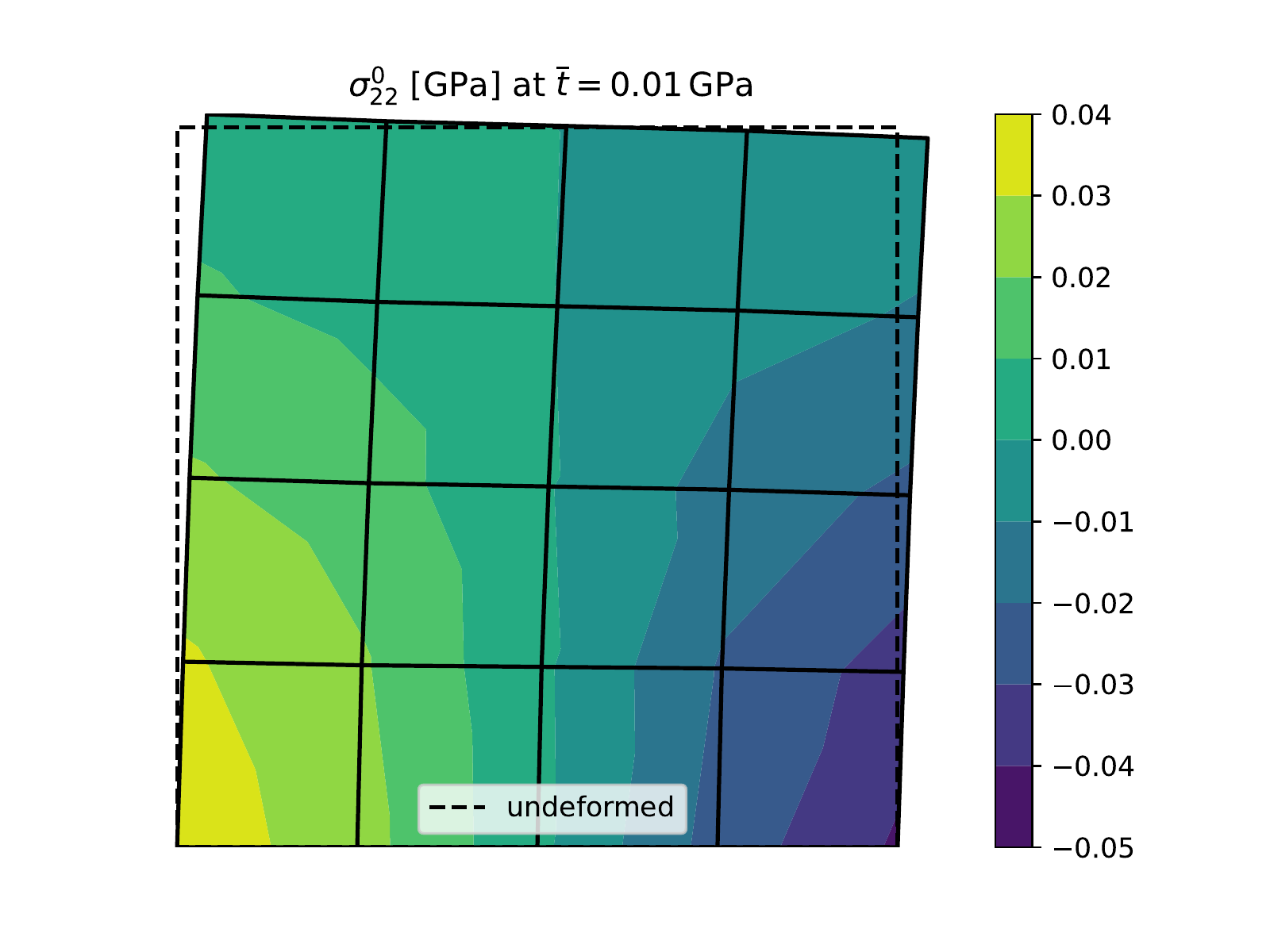} & 
  \includegraphics[width=.49\textwidth]{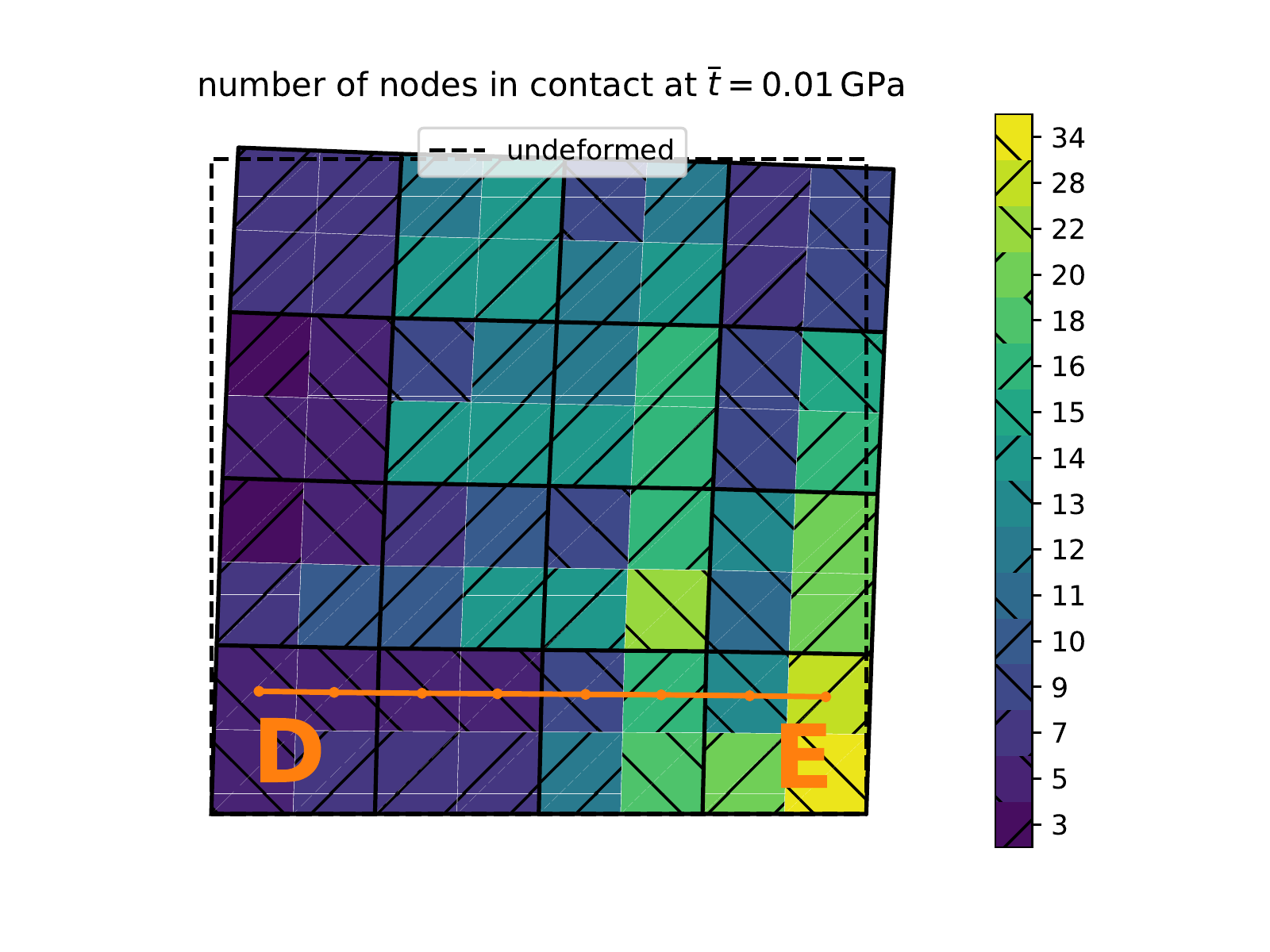}
  \end{tabular}
  \caption{Bending: Macroscopic results -- vertical stress component and number of nodes in contact.}
  \label{fig-bending-macro}
\end{figure}

\begin{figure}
  \centering
   \begin{tabular}{cc}
  \includegraphics[width=.49\textwidth]{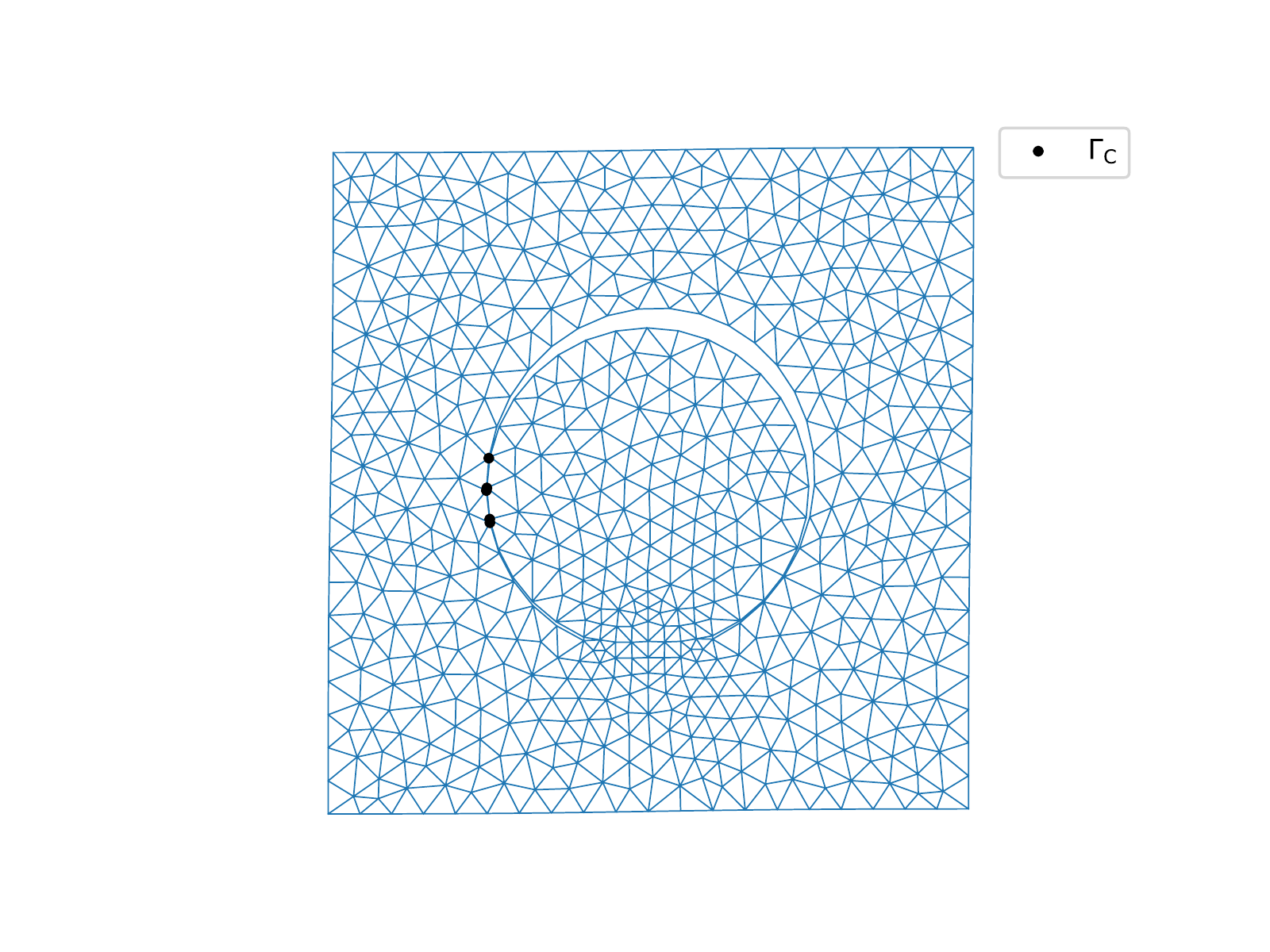}
  \includegraphics[width=.49\textwidth]{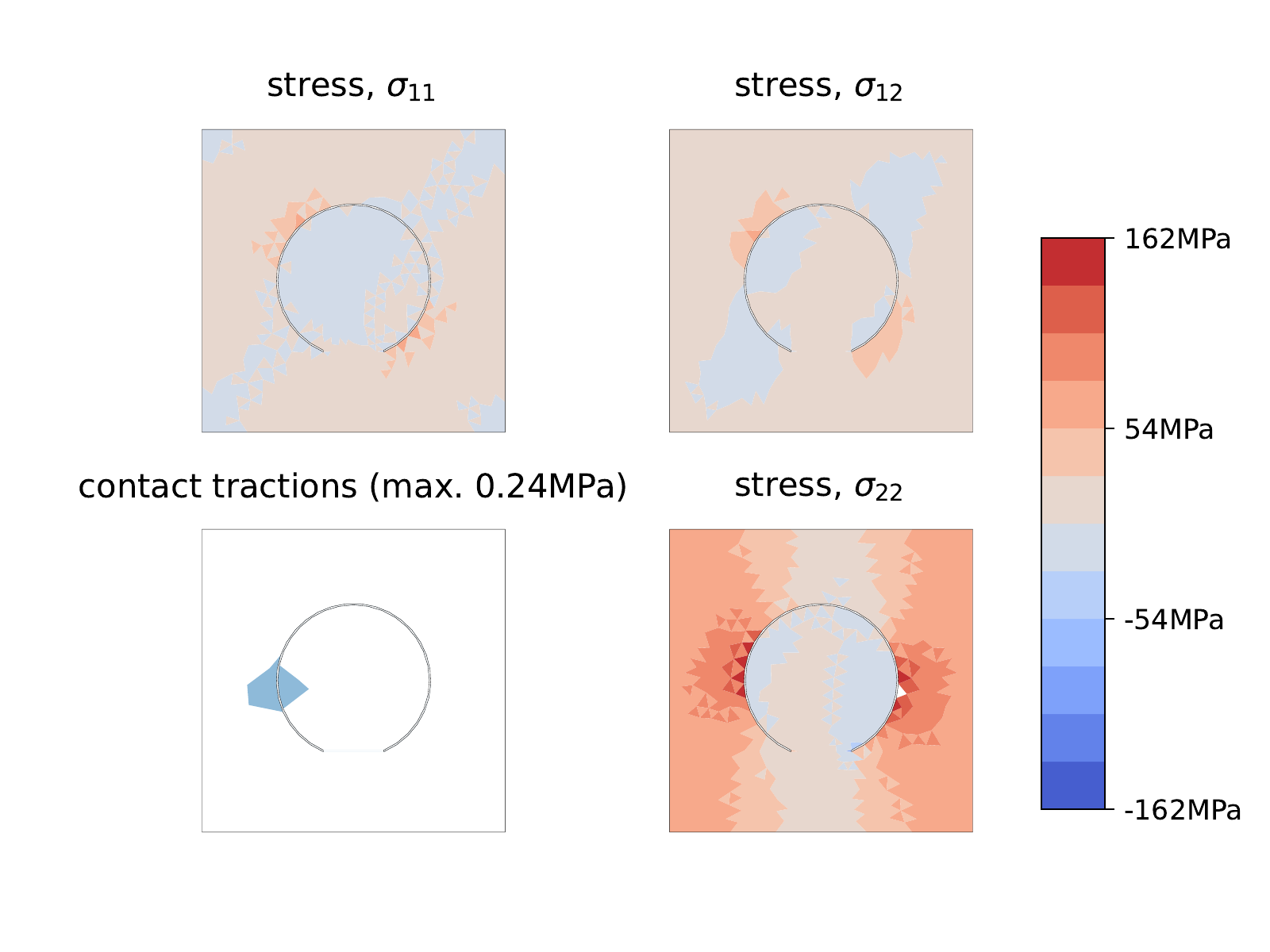}
  \end{tabular}
  \caption{Deformed shape and state of stress at the microscopic level (location D).}
  
  \label{fig-bending-micro-solution}
\end{figure}

Fig.~\ref{fig-bending-macro} shows the spatial distribution of the vertical stress component and contact at the macroscopic level.
Fig.~\ref{fig-bending-micro-solution} depicts a deformed microscopic cell and distribution of stress components.

\section{Conclusion}\label{sec-conclusion}
We proposed and tested the new two-scale method called the Macroscopic Contact (MC) method for numerical modelling of the porous structures featured by local self-contact on the pore surfaces. Within the iterations associated with displacement increments, this method consists of the micro- and macro-level steps; the latter one is formulated as the contact problem with the ``two-scale'' contact surface $\Sigma_\Gamma$ constituted by a neighborhood of the local active contact surface, as defined in \eq{eq-ma-2*}. The MC method appears to converge faster than the method involving the macroscopic step in the form of linear elasticity problem with a consistent incremental modulus, as proposed in  \cite{Rohan-Heczko-CaS}. A combination of the two algorithms will be explored in a further research. Issues of the computation stability will require further effort especially when considering 3D structures and dynamic loading; in this respect, aproaches reported in \cite{Kopacka2018b,Kolman2021} will be followed.

\paragraph{Acknowledgment}
The research has been supported by the grant project GA~22-00863K of the Czech Science Foundation.

\begin{appendices}

\end{appendices}

\bibliographystyle{plain}      

\bibliography{biblio-homog-contact}

\end{document}